\documentclass{amsart}

\usepackage{amsmath,amsthm}
\usepackage{amsfonts,amssymb}

\newtheorem{thm}{Theorem}[section]
\newtheorem{cor}[thm]{Corollary}
\newtheorem{lem}[thm]{Lemma}
\newtheorem{prop}[thm]{Proposition}

\theoremstyle{remark}

 \def\xb{{\mathbf x}}

 \def\CA{{\mathcal A}}     
 \def\CB{{\mathcal B}}

 \def\CH{{\mathcal H}}
 \def\CI{{\mathcal I}}

 \def\CR{{\mathcal R}}

 \def\CV{{\mathcal V}}

 \def\NN{{\mathbb N}}

 \def\RR{{\mathbb R}}

 \def\proj{\operatorname{proj}}

\begin{document}
\title[]
{A Discretized Fourier Orthogonal Expansion in Orthogonal Polynomials on a Cylinder} 
 
\author{Jeremy Wade}
\address{Department of Mathematics\\224 Yates\\ Pittsburg State University\\1701 South Broadway\\
    Pittsburg, Kansas 66762.}
\email{jwade@pittstate.edu}

\date{\today}
\thanks{}

\begin{abstract}
	We study the convergence of a discretized Fourier orthogonal expansion in orthogonal polynomials on $B^2 \times [-1,1]$, where $B^2$ is the closed unit disk in $\RR^2$.  The discretized expansion uses a finite set of Radon projections and provides an algorithm for reconstructing three dimensional images in computed tomography.  The Lebesgue constant is shown to be $m \, (\log(m+1))^2$, and convergence is established for functions in $C^2(B^2 \times [-1,1])$.
\end{abstract}

\maketitle

\section{Introduction}
\setcounter{equation}{0}

The Radon projection of function on $\RR^d$, defined by
\begin{displaymath}
	\CR_\theta(f;t) = \int_{\langle x, \theta \rangle=t} \, f(x) \, dx
\end{displaymath}
for $t \in \RR$ and $\theta \in S^{d-1}$, the unit sphere in $\RR^d$, is a classical topic with many applications to approximation theory and image reconstruction.  In recent years, the role of Radon projections in topics dealing with orthogonal polynomials and Fourier orthogonal expansions in terms of orthogonal polynomials have been studied.  These results rely on results put forward in two papers in by Marr \cite{Marr} in 1974 and Logan and Shepp \cite{LS} in 1975, which uncovered a relationship between Radon projections and orthogonal polynomials on $B^2$, the unit disk in $\RR^2$.  This relationship has led to many results in multi-dimensional approximation theory involving Radon projections.  In 2005, Bojanov and Xu used this relationship to find a polynomial on $B^2$ that interpolates the Radon projections of a function taken on sets of parallel lines in directions given by equidistant angles along the unit circle \cite{BojanovXu}, while in 2004, Bojanov studied the polynomial interpolating Radon projections on parallel lines in arbitrary directions \cite{Bojanov}.  In 1998, Petrushev extended the relationship in \cite{LS} and \cite{Marr} to a connection between ridge polynomials on $B^d$ and higher-dimensional Radon projections \cite{Petrushev}.  In 2006, Xu introduced a discretized Fourier orthogonal expansion on $B^2$, where the discrete data are finite Radon projections, and then studied the convergence properties of this expansion in the uniform norm \cite{XuRadon}.  This discretized expansion has application to image reconstruction in computerized tomography, where X-ray data corresponds to Radon projections.  This discretized Fourier orthogonal expansion was generalized to the domain $B^d$ in \cite{XuBd}.

In this paper, we study a discretization of the Fourier orthogonal expansion of functions on the cylinder $B^2 \times [-1,1]$ in terms of orthogonal polynomials with respect to the measure $w(x,y,z)=(1-z^2)^{-1/2}$, where $(x,y) \in B^2$ and $z \in [-1,1]$.  The discrete data are Radon transforms, taken on parallel disks, which are perpendicular to the axis of the cylinder.  This discretized expansion is a linear operator on $C(B^2 \times [-1,1])$, and the Lebesgue constant is shown to be $m(\log(m+1))^2$.  While the construction of this discretized expansion uses on the results from the discretized expansion on $B^2$ in \cite{XuRadon}, the proof of the Lebesgue constant of the discretized expansion utilizes a different approach.  This is due to our definition of the degree of a polynomial being the \emph{total degree}, which precludes us from expressing the orthogonal polynomial series on $B^2 \times [-1,1]$ as a simple product of two orthogonal polynomial series on $B^2$ and $[-1,1]$, and hence from obtaining a result trivially from previous results on these regions.  In estimating the upper bound of the Lebesgue constant, generating functions for orthogonal polynomial series are used to separate the estimates on $B^2$ and $[-1,1]$.  This approach also provides an alternative proof to the original results on $B^2$.  The discretized algorithm has application in 3-d computerized tomography, as the discrete Radon projections can be interpreted as discrete X-ray data.

This paper is organized as follows.  In section 2, background material is presented, including some information on orthogonal polynomials and the Radon projection.  The Fourier orthogonal expansion is briefly discussed in section 3, and the algorithm is presented in section 4.  Since the derivation of the algorithm is also contained in \cite{XuRadon}, sections 3 and 4 are brief.  A detailed proof of the Lebesgue constant is contained in section 5.  

\section{Preliminaries}
\setcounter{equation}{0}
\subsection{Radon projections}
Let $f$ be a function defined on $\RR^2$.  We define $L(\theta,t)$, for $\theta \in [0,2\pi]$ and $t \in \RR$, to be the line $\{ (s \cos \theta + t \sin \theta, s \sin \theta - t \cos \theta)\ : \ s \in \RR \}$.  The Radon projection of $f$, $\CR_\theta(f;t)$ is defined to be
\begin{displaymath}
	\CR_\theta(f;t):=\int_{L(\theta;t)} \, f(x,y) \, dx dy.
\end{displaymath}

By restricting the domain of $f$ to $B^2$, we only need to consider line segments in $B^2$ in the definition of $\CR_\theta(f;t)$.  For this reason, we define $I(\theta,t)$ to be the intersection of $L(\theta;t)$ and $B^2$; that is,
\begin{displaymath}
	I(\theta;t)= \{ (s \cos \theta + t \sin \theta, s \sin \theta - t \cos \theta)\ : \ |s| \leq \sqrt{1-t^2}\ \},
\end{displaymath}
and re-define $\CR_\theta(f;t)$ with the integral taken over $I(\theta;t)$ instead of $L(\theta;t)$.

For functions on $B^2 \times [-1,1]$, we will use Radon projections on parallel disks, perpendicular to the axis of the cylinder.  For these functions, we define the Radon projection by
\begin{displaymath}
	\CR_\theta(f (\cdot, \cdot, z);t) := \int_{I(\theta;t)}\, f(x,y,z) \, dx dy,
\end{displaymath}
where $(x,y) \in B^2$ and $z \in [-1,1]$.

Radon's famous result states that a function can be completely reconstructed from its Radon projections, provided the function satisfies some modest decay conditions.  Radon gave an explicit reconstruction formula, which required all the Radon projections of a function to be known.  We instead seek to approximate a function using only a finite number of Radon projections.

\subsection{Orthogonal polynomials}
Let $\Omega$ be a subset of $\RR^d$ with positive Lebesgue measure, and $\xb = (x_1, x_2, \ldots, x_d) \in \Omega$.  Let $\alpha = (\alpha_1, \alpha_2, \ldots, \alpha_d) \in \NN_0^d$, and define
\begin{displaymath}
	\xb^\alpha = x_1^{\alpha_1} x_2^{\alpha_2} \cdots x_d^{\alpha_d}, \ \ |\alpha| = \alpha_1 + \alpha_2 + \cdots + \alpha_d.
\end{displaymath}
We define a polynomial $P$ in $d$ variables to be of degree $n$ if it is of \emph{total degree} $n$, which means $P$ is of the form
\begin{displaymath}
	P(\xb) = \sum_{|\alpha| \leq n} c_\alpha \xb^\alpha,
\end{displaymath}
where $c_\alpha$ are real numbers, and at least one $c_\alpha$ with $|\alpha|=n$ is not zero.

Let $\omega$ be a positive measure on $\Omega$ with finite moments; that is, all integrals of monomials are finite.  A polynomial $P$ of degree $n$ is orthogonal on $\Omega$ with respect to $\omega$ if
\begin{displaymath}
	\int_\Omega \, P(x) Q(x) \, d\omega(x) = 0
\end{displaymath}
whenever $Q$ is a polynomial of degree less than $n$.

We denote by $\CV_n(\Omega; \omega)$ the space of orthogonal polynomials of degree $n$ on $\Omega$ with respect to the weight function $\omega$.  In the case that $\omega = 1$, we simply write $\CV_n(\Omega)$.  It is well known that $\text{dim}(\CV_n(\Omega)) = \binom{n+d-1}{d}$; see, for example, \cite{DX}.  We define $\Pi_n$ to be space of polynomials of degree less than or equal to $n$.  We may write $\Pi_n = \sum_{k=0}^n \CV_k(\Omega)$, and it follows that $\text{dim } \Pi_n = \binom{n+d}{d}$.  An \emph{orthogonal polynomial sequence} is an ordered basis of orthogonal polynomials for the space of polynomials.   

On the domain $\Omega=[-1,1]$, the Chebyshev polynomials of the first and second kinds are denoted by
\begin{displaymath}
	T_n(x) = \cos n\theta,\ \ \ U_n(x) = \frac{\sin\left( (n+1)\theta\right)}{\sin \theta},
\end{displaymath}
respectively, with $x = \cos \theta$.  The polynomials $T_n(x)$ are orthogonal with respect to the weight function $\omega(x)=\tfrac{1}{\pi}(1-x^2)^{-1/2}$, and the polynomials $U_n(x)$ are orthogonal with respect to the weight functions $\tfrac{1}{\pi}(1-x^2)^{1/2}$.  We denote by $\tilde{T}_n$ the orthonormal Chebyshev polynomials of the first kind.  They satisfy
\begin{displaymath}
	\widetilde{T}_0(x) = T_0(x), \ \ \widetilde{T}_n(x) = \sqrt{2}T_n(x).
\end{displaymath}
The zeroes of $T_n(x)$ are 
\begin{equation}
	z_{l,n}:=\cos(\gamma_{l,n}) := \cos \left( \frac{2l+1}{2n}\pi \right), \ \ l=0,1,\ldots n-1,
	\label{cheb1zeroes}
\end{equation}
and the $n$-point Gaussian quadrature associated with the Chebyshev polynomials of the first kind is given by
\begin{equation}
	\frac{1}{\pi}\int_{-1}^1 \, f(x) \frac{dx}{\sqrt{1-x^2}} 
	\approx \frac{1}{n} \sum_{l=0}^{n-1} f \left( z_{l,n} \right).
	\label{quadraturecheb1}
\end{equation} 
It is well known that equality is obtained in \eqref{quadraturecheb1} if the degree of $f$ is less than or equal to $2n-1$.

The zeroes of $U_n(x)$ are
\begin{equation}
	\cos \theta_{j,n} := \cos \frac{j\pi}{n+1},\ \ j=1,2,\ldots n,
\label{cheb2zeroes}
\end{equation}
and the Gaussian quadrature associated with the Chebyshev polynomials of the second kind is given by
\begin{equation}
	\frac{1}{\pi} \int_{-1}^1 \, f(x) \sqrt{1-x^2} \, dx 
	\approx \frac{1}{n} \sum_{j=1}^{n} f \left( \cos\left(\theta_{j,n}\right)\right).
\label{quadraturecheb2}
\end{equation}
Again, equality is achieved in \eqref{quadraturecheb2} if $f$ is a polynomial of degree $2n-1$ or less. 

For $\Omega = B^2$,  it was shown in \cite{LS} that the polynomials
\begin{equation}
	U_n(\theta_{j,n}; x, y) := U_n(x\cos \theta_{j,n} + y \sin \theta_{j,n}),\ \ \ j=0,1,\ldots, n
	\label{B2basis}
\end{equation}
form an orthonormal basis for $\CV_n(B^2, \omega)$ with respect to the measure $\omega(x,y)=\frac{1}{\pi}$.  This basis plays an important role in the construction of the discretized Fourier expansion on $B^2$.  The remarkable relationship between Radon projections and orthogonal polynomials, discovered by Marr in 1974, is given below.
\begin{thm}
	\label{polypresB2}
	\cite[Thm. 1]{Marr}
	If $P \in \CV_m(B^2)$, then
	\begin{displaymath}
		\CR_\theta(P;t) =  2\frac{\sqrt{1-t^2}}{m+1} U_m(t) P(\cos \theta, \sin \theta),
	\end{displaymath}
	where $U_m(t)$ is the Chebyshev polynomial of the second kind of degree $m$.
\end{thm}

On $B^2 \times [-1,1]$, define the weight function $\omega_Z(x,y,z) := \tfrac{1}{\pi}(1-z^2)^{-1/2}$, with $(x,y) \in B^2$ and $z \in [-1, 1]$.  It follows from the above facts that, for $(x,y) \in B^2$ and $z \in [-1,1]$, the polynomials
\begin{align}
	& P_{j,k,n}(x,y,z):= U_k(\theta_{j,k};x,y)\widetilde{T}_{n-k}(z) \ , \ \ k=0,1,\ldots n,\ j=0,1,\ldots k \label{onpolycyl}
\end{align}
form an orthonormal basis for $\CV_n(B^2 \times [-1,1];\omega_Z)$.  This basis plays an important role in the construction of the discretized Fourier orthogonal expansion on $B^2$ and $B^2 \times [-1,1]$.

\section{Fourier Orthogonal Expansion and Main Result}
\setcounter{equation}{0}

Standard Hilbert space theory allows us to decompose $L^2(\Omega;\omega)$ as
\begin{displaymath}
	L^2(\Omega;\omega) = \bigoplus_{k=0}^\infty \CV_k(\Omega;\omega), \ \ f = \sum_{k=0}^\infty \proj_k f, 
\end{displaymath}
provided $\omega$ is a positive measure with finite moments on $\Omega$.  The Fourier partial sum of degree $m$, $S_n$, is defined by
\begin{displaymath}
	S_n f = \sum_{k=0}^n \, \proj_k f.
\end{displaymath}
It is a standard result that $S_n f$ is the best approximation to $f$ in $\Pi_n$ in the norm of $L^2(\Omega;\omega)$.  We will use these standard results to construct our discretized expansion on the cylinder below.  We first list the results on $B^2$, and then move on to the cylinder.

\subsection{Discretized Expansion on $B^2$}

The results in this sub-section are due to Xu and are proven in \cite{XuRadon}.  Recall the orthogonal polynomials $U_k(\theta_{j,k};x,y)$ in \eqref{B2basis} form an orthonormal basis for $\CV_k(B^2)$.  Using Marr's theorem \cite{Marr}, Xu was able to express the Fourier coefficients of the Fourier orthogonal expansion of $f$ in this basis in terms of the Radon projections of $f$.
\begin{thm}
	\label{RadonCheb}
	\cite[Prop. 3.1]{XuRadon}
	Let $m >0$ and $f \in L^2(B^2)$.  For $0 \leq k \leq 2m$ and $0 \leq j \leq k$, 
	\begin{multline*}
		\int_{B^2} f(x,y) U_k(\theta_{j,k};x,y) \, dx dy \\
		= \frac{1}{2m+1}\sum_{\nu=0}^{2m} \frac{1}{\pi} \int_{-1}^1 \, \CR_{\phi_\nu}(f;t), U_k(t) \, dt \, U_k(\cos(\theta_{j,k}-\phi_\nu)),
	\end{multline*}
	where $\phi_\nu = \frac{\nu}{2m+1}$.
\end{thm}
As a result of this relationship, the orthogonal projection of a function $f$ in $L^2(B^2)$ onto $\CV_k(B^2)$ can be written in terms of the Radon projections of $f$.
\begin{thm}
	\label{B2proj}
	\cite[Thm. 3.2]{XuRadon}
	For $m >0$ and $k \leq 2m$, the projection operator $\proj_k$ from $L^2(B^2)$ to $\CV_k(B^2)$ can be written in the form
	\begin{displaymath}
		\proj_k f (x,y) = \frac{1}{2m+1} \sum_{\nu=0}^{2m} \frac{1}{\pi} \int_{-1}^1 \, \CR_{\phi_\nu}(f;t) U_k(t) \, dt \, (k+1) U_k(\phi_\nu;x,y).
	\end{displaymath}
\end{thm}
As a consequence of this theorem, the Fourier partial sum of $f$, $S_{2m} f$, can also be written in terms of the Radon projections of $f$,
\begin{equation}
	\label{B2FourierPartial}
	S_{2m} f (x,y) = \frac{1}{2m+1} \sum_{k=0}^{2m} \sum_{\nu=0}^{2m} \frac{1}{\pi} \int_{-1}^1 \, \CR_{\phi_\nu}(f;t) U_k(t) \, dt \, (k+1) U_k(\phi_\nu;x,y).
\end{equation}
Recalling \eqref{polypresB2}, the expression
\begin{displaymath}
	\frac{\CR_{\phi_\nu}(f;t)}{\sqrt{1-t^2}}
\end{displaymath}
is a polynomial of degree $2m$ if $f$ is a polynomial of degree $2m$.  Hence, by multiplying and dividing by a factor of $\sqrt{1-t^2}$ in \eqref{B2FourierPartial}, and then using the $2m$-point Gaussian quadrature rule given in \eqref{quadraturecheb2} to replace the integral with a sum, a discretized Fourier orthogonal expansion which preserves polynomials of degree less than or equal to $2m-1$ is obtained.  This discretized expansion, $\CA_{2m}$, is given by
\begin{align*}
	\CA_{2m}(f)(x,y) = & \frac{1}{(2m+1)^2} \sum_{\nu=0}^{2m}  \sum_{j=1}^{2m} 
	\CR_{\phi_\nu}(f;\cos \theta_{j,2m}) \\
	& \qquad \times \sum_{k=0}^{2m} (k+1) \sin( (k+1) \theta_{j,2m} ) U_k(\phi_\nu;x,y).
\end{align*}
As an operator on $C(B^2)$, $\CA_{2m}$ has a Lebesgue constant of $m\log(m+1)$.  The upper limit of the sum in $\nu$ is chosen to be $2m$ to eliminate redundancy in the Radon data.  With the choice of $2m$, the Radon projections are taken along parallel lines in $(2m+1)$ directions, given by equally spaced points on the unit circle.  If the Radon projections were taken along parallel lines in $2m$ directions, then for $\nu < m$, and $\phi_\nu = \frac{2\pi\nu}{2m}$,
\begin{displaymath}
	\pi + \phi_nu= \phi_{\nu+m}, 
\end{displaymath}
and the identity
\begin{displaymath}
	\CR_{\pi+\phi_\nu}(f;\cos(\theta_{2m+1-j,2m})) = \CR_{\phi_\nu}(f;\cos(\theta_{j,2m}))
\end{displaymath}
shows these two Radon projections are the same.  In practical settings, more Radon projections are desirable, so we choose to take Radon projections in an odd number of directions.

\subsection{Discretized Expansion on $B^2 \times [-1,1]$}

The construction discretized Fourier orthogonal expansion uses the basis for $\CV_n(B^2 \times [-1,1]; \omega_Z)$ in \eqref{onpolycyl} and the results on $B^2$.

\begin{thm}
	Let $m \geq 0$ and let $n \leq 2m$.  Define $\phi_\nu := \frac{2\nu\pi}{2m+1}$, and $\sigma_\nu(x,y) := \arccos (x \cos(\phi_\nu) + y \sin(\phi_\nu))$.  The operator $\proj_n$ can be written as
	\begin{align}
		&\proj_n f(x,y,z) = 
		\frac{1}{\pi} \frac{1}{2m+1}
		\sum_{\nu=0}^{2m} 
	 	\int_{-1}^1 \int_{-1}^1 \CR_{\phi_\nu}(f(\cdot, \cdot, s);t)
		\Psi_\nu(x,y,z;s,t) \, dt \,\frac{ds}{\sqrt{1-s^2}},\\
		&\Psi_{\nu,n} (x,y,z;s,t)=\frac{1}{2m+1} \sum_{k=0}^n (k+1)
		U_k(t) U_k(\cos(\sigma_\nu(x,y))) \tilde{T}_{n-k}(s)\tilde{T}_{n-k}(z) \notag
	\label{cylproj}
	\end{align}
\end{thm}
\begin{proof}
	Since the polynomials $P_{j,k,n}$ in \eqref{onpolycyl} form an orthonormal basis for $\CV_n(B^2 \times [-1,1]; \omega_Z)$, 
	\begin{displaymath}
		\proj_n f(x,y,z) = \sum_{k=0}^n \sum_{j=0}^k \frac{1}{\pi} \int_{-1}^1 \int_{B^2} f(x,y,z) P_{j,k,n}(x,y,z) \omega_Z(x,y,z).
	\end{displaymath}
	After expanding $P_{j,k,n}(x,y,z)$, Theorem \eqref{B2proj} gives the result. 
\end{proof}

This relationship between the projection operator and the Radon projection also yields a connection between the partial sum operator $S_{2m}$ and the Radon projection.

\begin{cor}
	Let $m \geq 0$.  The partial sum operator $S_{2m}$ may be written as
	\begin{equation}
		S_{2m} (f) (x,y,z) = \sum_{\nu=0}^{2m}	\frac{1}{\pi} \int_{-1}^1 \int_{-1}^1
		\CR_{\phi_\nu} (f(\cdot, \cdot, s) ; t) \Phi_\nu(x,y,z;s,t) \, \frac{ds}{\sqrt{1-s^2}} \, dt ,
		\label{partialsumcyl}
	\end{equation}
where
	\begin{equation*}
		\Phi_\nu(x,y,z;s,t) = \sum_{n=0}^{2m}  \Psi_{\nu,n}(x,y,z;s,t).
	\end{equation*}
\end{cor}

We will discretize $S_{2m}$ by approximating the integrals with Gaussian quadratures.  In order to achieve a result on the convergence of our discretized expansion, we must first prove a result concerning the Radon projection of a polynomial on $B^2 \times [-1,1]$.
\begin{lem}
	\label{Radonpoly}
	If $P$ is a polynomial of degree $k$ on $B^2 \times [-1,1]$, then for $\theta \in [0, 2\pi]$,
		\begin{displaymath}
			\frac{ \CR_{\theta} (P (\cdot , \cdot, s); t)}{\sqrt{1-t^2}} 
		\end{displaymath}
	is a polynomial of degree $k$ in $t$.
\end{lem}
\begin{proof}
	If $P$ is a polynomial of degree $k$, we may write
	\begin{displaymath}
		P(x,y,z) = \sum_{i=0}^k c_i z^i p_{k-i}(x,y),
	\end{displaymath}
	where $p_{k-i}(x,y)$ is a polynomial of degree $k-i$ in $(x,y)$.  Following the proof of \cite[Lem. 2.2]{XuRadon}, we write
	\begin{displaymath}
		\frac{\CR_\theta( P(\cdot, \cdot, z);t) }{\sqrt{1-t^2}}:= \sum_{i=0}^k c_i z^i \frac{1}{\sqrt{1-t^2}} \int_{I(\theta;t)} \, p_{k-i}(x,y) \, dx dy.
	\end{displaymath}
	Rewriting the integral and changing variables,
	\begin{align*}
		& \frac{1}{\sqrt{1-t^2}} \int_{I(\theta;t)} p_{k-i}(x,y) \, dx dy\\ 
		& \qquad =\frac{1}{\sqrt{1-t^2}} \int_{-\sqrt{1-t^2}}^{\sqrt{1-t^2}} \, p_{k-i}(t \cos \theta + s \sin \theta, t \sin \theta - s \cos \theta )\, ds\\
		& \qquad = \int_{-1}^1 \, p_{k-i}(t \cos \theta + \sqrt{1-t^2} s \sin \theta, t \sin \theta - \sqrt{1-t^2} s \cos \theta )\, ds.
	\end{align*}
	After expanding $p_{k-i}(x,y)$ in the integrand, we note that each odd power of $\sqrt{1-t^2}$ is accompanied by an odd power of $s$, which becomes $0$ after integrating.  Hence, 
	\begin{align*}
		& \frac{1}{\sqrt{1-t^2}} \int_{I(\theta,t)} p_{k-i}(x,y) \, dx dy = \sum_{j=0}^{\lfloor \frac{k-i}{2} \rfloor} b_j t^{k-i-2j}(1-t^2)^j
&	\end{align*}
	for some coefficients $b_j$.  The lemma follows.
\end{proof}
With this lemma in mind, we substitue 
\begin{displaymath}
	\CR_{\phi_\nu}(f(\cdot, \cdot, s);t) = \frac{\CR_{\phi_\nu}(f(\cdot, \cdot, s);t)}{\sqrt{1-t^2}} \sqrt{1-t^2} 
\end{displaymath}
into \eqref{partialsumcyl}.  We then discretize the two integrals by using $2m$-point Guassian quadratures.  For the integral in $s$, the quadrature formula \eqref{quadraturecheb1} is used, and for the integral in $t$, we use the quadrature formula \eqref{quadraturecheb2}.  The discretized Fourier expansion, $\CB_{2m}$, is given below.
\begin{thm}  
	For $m \geq 0$, $(x,y) \in B^2$ and $z \in [-1,1]$, we define
	\begin{equation}
		\CB_{2m}(f)(x,y,z) :=  \sum_{\nu=0}^{2m} \sum_{j=1}^{2m} \sum_{l=0}^{2m-1}
		\CR_{\phi_\nu} (f(\cdot, \cdot, z_{l,2m}; \cos(\theta_{j,2m})) T_{\nu,j,l}(x,y,z),
		\label{B2m}
	\end{equation}
	where
	\begin{align}
		\label{T1}
		T_{\nu,j,l}(x,y,z) = & \frac{1}{(2m+1)^3} \sum_{n=0}^{2m} \sum_{k=0}^{n}(k+1) \sin\left(  (k+1) \theta_{j,2m} \right) \\
		\notag & \times U_k \left(\cos(\sigma_\nu(x,y))\right) T_{n-k}(z_l)T_{n-k}(z).
	\end{align}
\end{thm}

\subsection{Main Theorems}

As a result of lemma \ref{Radonpoly} and the fact that $2m$-point Gaussian quadratures are exact for polynomials of degrees up to $4m-1$, we obtain the following theorem.

\begin{thm}
	\label{polypres}
	The operator $\CB_{2m}$ preserves polynomials of degree less than or equal to $2m-1$; that is, for $(x,y) \in B^2$ and $z \in [-1,1]$,
	\begin{displaymath}
		\CB_{2m} (f) (x,y,z) = f(x,y,z)
	\end{displaymath}
	for $f \in \Pi_{2m-1}$.
\end{thm}

The proof of the following theorem is contained in the next section, and is the main substance of the paper.  Comparing this with the result in \cite{XuRadon}, we see that the extension of the algorithm to the cylinder introduces a factor of $\log (m+1)$.

\begin{thm}
	\label{Lebesgueconstant}
	For $m \geq 0$, the norm of the operator $\CB_{2m}$ on $C(B^2 \times [-1,1])$ is given by
	\begin{displaymath}
		\| B_{2m} \|_\infty \approx m \left( \log(m+1) \right)^2.
	\end{displaymath}
\end{thm}
The proof of this theorem is not trivial.  Since we have defined the degree of a polynomial to be its total degree, the series in the definition of $T_{\nu,j,l}(x,y,z)$,
\begin{displaymath}
	\sum_{n=0}^{2m} \sum_{k=0}^{n}(k+1) \sin\left(  (k+1) \theta_{j,2m} \right) U_k \left(\cos(\sigma_\nu(x,y))\right)
		\times \tilde{T}_{n-k}(z_l) \tilde{T}_{n-k}(z), 
\end{displaymath}
cannot be written as the product of a two series, one in terms of $z$ and $z_l$, and one in terms of $\theta_j$ and $\sigma_\nu(x,y)$.  As a result, the estimate of the Lebesgue constant cannot be trivially reduced to an estimate on $B^2$ and an estimate on $[-1,1]$.  In particular, for the upper bound of the estimate, a different approach from that which was used in \cite{XuRadon} is used to obtain our result. 

As a corollary of Theorems \eqref{Lebesgueconstant} and \eqref{polypres}, we obtain the following corollary.

\begin{cor}
	For $f \in C^2(B^2 \times [-1,1])$, $\CB_{2m} (f)$ converges to $f$ in the uniform norm.
\end{cor}
\begin{proof}
	If $f \in C^2(B^2 \times [-1,1])$, then by Theorem 1 in \cite{Bagby}, there exists a polynomial $p_n$ of degree $n$ on $B^2 \times [-1,1]$, and a constant $C>0$, so that
	\begin{displaymath}
		\|f-p_n\|_\infty \leq \frac{C}{n^2} \omega_{f,2} \left( \frac{1}{n}\right),
	\end{displaymath}
	where
	\begin{displaymath}
		\omega_{f,2} \left( \frac{1}{n} \right) = \sup_{|\gamma|=2} 
		\left( \sup_{\substack{  x,y \in B^2 \times [-1,1]\\
					|x-y|\leq 1/n}}
		\left| D^\gamma f(x)-D^\gamma(y)\right|\right).
	\end{displaymath}
	We let $n = 2m-1$ to obtain
	\begin{align*}
		\|\CB_{2m}(f) - f \|_\infty & \leq \| \CB_{2m} (f-p_{2m-1}) \|_\infty + \| f- p_{2m-1} \|_\infty\\
		& \leq \|f -p_{2m-1}\|_\infty \left( 1 + \|\CB_{2m}\|_\infty \right)\\
		& \leq c \, \frac{1}{(2m-1)^2} \left( m(\log(m+1))^2 + 1\right),
	\end{align*}
	which converges to zero as $m$ approaches infinity.
\end{proof}
Before proceeding to the proof of Theorem \eqref{Lebesgueconstant}, we make one comment.  We believe that the Lebesgue constant of the Fourier partial sum of the orthogonal expansion, $\|S_{2m}\|$, is $\approx m \log(m+1)$.  We are able to prove $\|S_{2m}\|_\infty = O(m\log(m+1))$, but we have yet to prove the lower bound.  If our intuition is correct, the discretization of the expansion adds a factor of $\log(m+1)$.  

\section{Proof of Theorem \ref{Lebesgueconstant}}
\setcounter{equation}{0}
We first derive an expression for which we may estimate $\|\CB_{2m}\|_\infty$.
\begin{prop}
	The norm of $\mathcal{B}_{2m}$ as an operator on $C(B^2\times [-1,1])$ is given by
\begin{equation}\label{BNorm}
	\|\mathcal{B}_{2m}\|_\infty=2 \max \sum_{\nu=0}^{2m} \sum_{j=1}^{2m}\sum_{l=0}^{2m-1} \sin \theta_{j,2m}  
	\left|T_{\nu,j,l}(x,y,z)\right| 
\end{equation}
where the maximum is taken over all points $(x,y,z)$ in $B^2\times[-1,1].$
\end{prop}\begin{proof}
	By definition, 
\begin{multline}
	\mathcal{R}_{\phi_\nu} ( f( \cdot, \cdot, z_l ) , \cos \theta_{j,2m} )
	= \int_{ I ( \cos \theta_{j,2m} , \phi_\nu ) } f(\tilde{x},\tilde{y},z_l)\, d\tilde{x}d\tilde{y}\\
	= \int_{ -\sin \theta_{j,2m} }^{ \sin \theta_{j,2m} } \, 
	f( \cos \theta_{j,2m} \cos \phi_\nu -s \cos \phi_\nu , \cos \theta_{j,2m} \sin \phi_\nu +s \cos \phi_\nu , z_l )
	\, ds.
\end{multline}	
Taking absolute value of both sides and using the triangle inequality, we immediately have
\begin{displaymath}
	\|\mathcal{B}_{2m}\|_\infty
	\leq 2 \max \sum_{\nu=0}^{2m} \sum_{j=1}^{2m} \sum_{l=0}^{2m-1} 
	\sin \theta_{j,2m} \left|T_{\nu,j,l}(x,y,z)\right| 
\end{displaymath}
On the other hand, if we define
\begin{displaymath}
	\mathcal{T}(x,y,z) := 2 \sum_{\nu=0}^{2m} \sum_{j=1}^{2m}\sum_{l=0}^{2m-1} 
	\sin \theta_{j,2m} \left|T_{\nu,j,l}(x,y,z)\right|, 
\end{displaymath}
then $\mathcal{T}(x,y,z)$ is a continuous function on $B^2 \times [-1,1]$, and hence achieves its maximum at some point $(x_0,y_0,z_0)$ on the cylinder.  We would like to choose a function $f$ so that $f(x,y,z)=\text{sign}(T_{\nu,j,l}(x_0,y_0,z_0))$ on the set of lines $\{(I( \cos \theta_{j,2m} ,\phi_\nu), z_l)\}$, for $1 \leq j \leq 2m$, $0 \leq \nu \leq 2m$, and $0 \leq l \leq 2m-1$, since this would immediately give us the result.  However, such a function may not be continuous at the points of intersection of these lines.  To allow for continuity, we instead take neighborhoods of volume $\varepsilon$  around each point of intersection of the lines, and define a function $f^*$ which is equal to $\text{sign}(T_{\nu,j,l}(x_0,y_0,z_0))$ on the lines $\{(I(\cos \theta_{j,2m} ,\phi_\nu),z_l) \}_{j,\nu,l}$ except on the $\varepsilon$-neighborhoods at the points of intersection; on the rest of the cylinder, $f^*$ is chosen so that it takes values between $1$ and $-1$ and is continuous.  It then follows that
\begin{displaymath}
	\left\|\mathcal{B}_{2m}\right\|_\infty
	\geq \left|\mathcal{B}_{2m}(f^*(x_0,y_0,z_0))\right| 
	\geq 2 \sum_{\nu=0}^{2m} \sum_{j=1}^{2m} \sum_{l=0}^{2m-1} 
	\sin \theta_{j,2m} \left|T_{\nu,j,l}(x_0,y_0,z_0)\right|-c\varepsilon,
\end{displaymath}
where $c$ denotes the number of points of intersection of the lines $\{( I(\cos \theta_{j,2m}, \phi_\nu), z_l)  \}_{j,\nu,l}$.  Since $\varepsilon$ is arbitrary, this proves the proposition.
\end{proof}

For the remainder of the proof, the number $n$ in \eqref{cheb1zeroes} and \eqref{cheb2zeroes} will be fixed as $2m$.  For this reason, we define 
\begin{align}
	\label{angledefns}
	&\theta_j = \theta_{j,2m} = \frac{j\pi}{2m+1},\  \gamma_l = \gamma_{l,2m} = \frac{2l+1}{4m}\pi, \  z_l=z_{l,2m}=\cos \gamma_l, \\
	& \phi_\nu = \frac{2\pi\nu}{2m+1}, \ \sigma_\nu(x,y) = \arccos( x\cos \phi_\nu + y \sin \phi_\nu) \notag.
\end{align}

The proof will be separated into two parts, a lower and an upper bound. 

\subsection{Lower Bound}

We will establish there exists a constant $c > 0$ so that $\|\CB_{2m} \|_\infty \geq c m (\log(m+1))^2$ for all $m > 0$.
By \eqref{BNorm}, it suffices to show 
\begin{displaymath}
	\sum_{\nu=0}^{2m} \sum_{j=1}^{2m} \sum_{l=0}^{2m-1} \sin \theta_j \left| T_{\nu,j,l}(x_1,y_1,z_1) \right| \geq c_1 m (\log(m+1))^2
\end{displaymath}
for the point $(x_1,y_1,z_1)=\left( \cos \tfrac{\pi}{4m+2}, \sin \tfrac{\pi}{4m+2}, 1 \right)$ for some $c > 0$.  We begin by deriving a compact formula for $T_{\nu,j,l}(x_1,y_1,1)$. Using the Christoffel-Darboux formula for $\widetilde{T}_n$, letting $\cos \gamma_z = z$,
\begin{align*}
	& \left| T_{\nu,j,l}  (x,y,z) \right| = \frac{1}{(2m+1)^3}\left| \frac{1}{\sin \sigma_\nu(x,y) } 
	\sum_{k=0}^{2m} \sin ( (k+1) \sigma_\nu (x,y) ) \sin( (k+1) \theta_{j} )  \right. \\
	& \left. \qquad \times \frac{ \cos ( (2m-k+1) \gamma_z ) \cos ( (2m-k) \gamma_l )
	- \cos ( (2m-k+1) \gamma_l ) \cos( (2m-k) \gamma_z )}{\cos ( \gamma_z ) 
	- \cos (\gamma_l). } \right|.
\end{align*}
Substituting in $z=1$ and applying the identity for the difference of cosines,
\begin{align*}
	& \left| T_{\nu,j,l}  (x,y,1) \right| 
	=  \frac{1}{(2m+1)^3}\frac{1}{\sin \sigma_\nu(x,y) }\frac{1}
	{\sin  \tfrac{\gamma_l}{2} } \\
	& \qquad \times \left| \sum_{k=0}^{2m} (k+1) 
	\sin ( (k+1) \sigma_\nu (x,y) ) \sin( (k+1) \theta_{j} )  
	\sin \left( (2m-k+1/2) \gamma_l \right) \right|.
\end{align*}
Applying the product formula for sine and the product formula for sine and cosine, 
\begin{align} 
	\label{T2} 
	&\left| T_{\nu,j,l} (x,y,1) \right| = 
	\frac{1}{4} \frac{1}{ (2m+1)^3 }
	\frac{1} { \sin \tfrac{ \gamma_l } { 2 } } 
	\frac{1} { \sin \sigma_\nu(x,y)  }\\
	& \notag \quad 
	\times \bigg| \sum_{k=0}^{2m} (k+1) \left[ 
	\sin ( (k+1)(\theta_{j}-\sigma_\nu(x,y) + \gamma_l) - \tfrac{3}{2}\gamma_l -\tfrac{\pi}{2} ) \right.\\
	& \notag \qquad \qquad \qquad \ 
	- \sin ( (k+1)(\theta_{j} - \sigma_\nu(x,y) - \gamma_l) + \tfrac{3}{2}\gamma_l + \tfrac{\pi}{2} ) \\ 
	& \notag \qquad \qquad \qquad \ 
	-\sin ( (k+1)(\theta_{j} + \sigma_\nu(x,y) + \gamma_l) - \tfrac{3}{2}\gamma_l -\tfrac{\pi}{2} ) \\
	& \notag \left. \qquad \qquad \qquad \ 
	+ \sin ( (k+1)(\theta_{j} + \sigma_\nu(x,y) - \gamma_l) + \tfrac{3}{2}\gamma_l + \tfrac{\pi}{2} ) 
	\right] \bigg|.\notag 
\end{align}
Next, apply the formula
\begin{align*}
	& \sum_{k=0}^{2m} (k+1) \sin( (k+1) \theta + \phi ) \\
	& \qquad \qquad = \frac{1}{2} \frac{ (2m +2) \sin( (2m+1) \theta + \phi )
	-(2m+1) \sin( (2m+2) \theta + \phi ) + \sin(\phi)}{\sin^2\left( \frac{\theta}{2} \right)}\\
	& \qquad \qquad = \frac{1}{2} \frac{\sin( (2m+1)\theta + \phi) -(4m+2)\cos( (2m+3/2)\theta + \phi)\sin(\theta/2)
	+ \sin(\phi)}{\sin^2\left( \frac{\theta}{2} \right)},
\end{align*}
 to \eqref{T2}.  Under our choice of $x_1$ and $y_1$, $\cos \sigma_\nu(x_1,y_1)  = \cos \tfrac{2\nu - 1/2}{2m+1}\pi $, so $\sigma_\nu(x_1,y_1) = \tfrac{2\nu-1/2}{2m+1}\pi$ if $1 \leq \nu \leq m$.  We will only be considering $\nu$ within this range, so we define $\sigma_\nu :=  \frac{2\nu-1/2}{2m+1}\pi$.  
Define
\begin{equation*}
	F_{j,l}^{\pm}(\theta,\phi,\gamma) 
	:= \frac{(-1)^{j+l+1}\cos\left(\frac{\gamma}{2} \right) 
	+(-1)^{j+l+1} (4m+2) \sin \left( \frac{\theta+\phi}{2} \right)
	  \sin \left(\frac{\theta + \phi + \gamma}{2} \right)
	  \pm \cos \left( \frac{3\gamma}{2} \right)}
	  { \sin^2 \left( \frac{\theta + \phi + \gamma }{2}\right)}.
\end{equation*}
Taking into account
\begin{align*}
	&(2m+1)(\theta_{j,2m} \pm_1 \sigma_\nu \pm_2 \gamma_l) = (j \pm_1( 2\nu - 1/2)\pm_2 (l + 1/2 ))\pi \pm_2 \gamma_l,
\end{align*}
where the subscripts indicate the signs of $\pm_1$ and $\pm_2$ are not related (a convention we will adopt for the remainder of the paper), we are able to write
\begin{align}\label{T3}
	&\left|T_{\nu,j,l} \left( \cos \tfrac{\pi}{4m+2} , \sin  \tfrac{\pi}{4m+2} ,1 \right)\right| = 
	\frac{1}{(2m+1)^3}\frac{1}{8\sin \sigma_\nu \sin \tfrac{\gamma_l}{2} }\\
	& \times \left|
	F_{j,l}^+(\theta_{j},-\sigma_\nu,\gamma_l)
	-F_{j,l}^-(\theta_{j},\sigma_\nu,\gamma_l)
	-F_{j,l}^-(\theta_{j},-\sigma_\nu,-\gamma_l) 
	+F_{j,l}^+(\theta_{j},\sigma_\nu,-\gamma_l) \right| \notag.  
\end{align} 
We will show the lower bound is attained if we restrict the summation in \eqref{BNorm} to the set of indices where $\pi/4 \leq \theta_{j} \leq 3\pi/8$, $\pi/4 \leq \gamma_l < \theta_{j}$, and $0 \leq \sigma_\nu < \theta_{j} - \gamma_l$, so we only take the sums over the following range of indices:
\begin{itemize}
	\item $\lfloor \frac{m}{2} \rfloor + 5 \leq j \leq 3 \lfloor \frac{m}{4} \rfloor$
	\item $\lfloor \frac{m}{2} \rfloor + 1 \leq l \leq j-4$
	\item $1 \leq \nu \leq \lfloor \frac{j-l}{2} \rfloor -1 $
\end{itemize}
We assume $m \geq 24$, so that these inequalities make sense.  With this restriction of summation, $\sin \theta_{j} $ and $\sin \tfrac{\gamma_l}{2}$ are bounded away from zero by a positive constant.  Hence, we are left with proving the estimate
\begin{align}
	\label{lowerbound}
	&\sum_{j=\lfloor \frac{m}{2} \rfloor +5}^{3\lfloor \frac{m}{4} \rfloor} 
	\sum_{l=\lfloor \frac{m}{2} \rfloor +1}^{j-4 } 
	\sum_{\nu=1}^{ \lfloor \frac{j-l}{2} \rfloor -1} 
	\frac{1}{(2m+1)^3}\frac{1}{\sin \sigma_\nu}\\
	\notag	& \qquad \times \left|
	F_{j,l}^+(\theta_{j},-\sigma_\nu,\gamma_l)
	-F_{j,l}^-(\theta_{j},\sigma_\nu,\gamma_l)
	-F_{j,l}^-(\theta_{j},-\sigma_\nu,-\gamma_l) 
	+F_{j,l}^+(\theta_{j},\sigma_\nu,-\gamma_l) \right|\\
	\notag &\qquad \qquad \qquad \geq c m (\log(m+1))^2 
\end{align}
Also note that, under this restriction of summation,
\begin{displaymath}
	0 < \frac{\theta_{j} -\sigma_\nu -\gamma_l}{2} \leq \frac{\theta_{j} + \sigma_\nu + \gamma_l}{2} \leq \frac{3}{8}\pi,
\end{displaymath}
so 
\begin{displaymath}
	\sin\left(\frac{\theta_{j} \pm_1 \sigma_\nu \pm_2 \gamma_l}{2}\right) \approx   \theta_{j} \pm_1 \sigma_\nu \pm_2 \gamma_l,
\end{displaymath}
where we have used the fact
\begin{displaymath}
	 \sin \theta \approx \theta
\end{displaymath}
if $-15\pi/16 \leq \theta \leq 15 \pi/16$, a fact we will use repeatedly throughout the proof.

The dominating terms in the summation will be the terms
\begin{displaymath}
	\frac{(4m+2)\sin \left(\tfrac{\theta_j \pm \sigma_\nu}{2}\right)}{\sin\left( \tfrac{\theta_j \pm \sigma_\nu -\gamma_l}{2} \right)},
\end{displaymath}
the middle term in the numerator of from $F_{j,l}^\pm(\theta_j, \pm \sigma_\nu, -\gamma_l)$.  We first prove two lemmas to eliminate the non-dominating terms.  The first lemma eliminates the first and third terms in the numerators of $F_{j,l}^{\pm_1} (\theta_j, \pm_{2}\sigma_\nu, \pm_3 \gamma_l)$.

\begin{lem}
	\label{lowerboundlem1}
	Recalling \eqref{angledefns},
	\begin{align*}
	J_1(\pm_1, \pm_2) &:= \frac{1}{(2m+1)^3} 
	\sum_{j=\lfloor \frac{m}{2} \rfloor +5}^{3\lfloor \frac{m}{4} \rfloor} 
	\sum_{l=\lfloor \frac{m}{2} \rfloor +1}^{j-4 } 
	\sum_{\nu=1}^{ \lfloor \frac{j-l}{2} \rfloor -1} 
	\frac{1}{\sin \sigma_\nu  \sin^2 \left( \tfrac{\theta_{j} \pm_1 \sigma_\nu \pm_2 \gamma_l}{2} \right)}	\\
	& \leq cm \log(m). 
\end{align*}
\end{lem}
\begin{proof}
	First, considering $\theta_{j} \pm \sigma_\nu + \gamma_l$, apply the inequalities
\begin{displaymath}
	\theta_{j}+ \sigma_\nu + \gamma_l > \theta_{j} + \gamma_l , \text{ and } \ \theta_{j} - \sigma_\nu + \gamma_l > \theta_{j} + \gamma_l - \pi/8,
\end{displaymath}
to obtain 
\begin{align*}
	J_1(\pm_1, +) 
	& \leq \frac{1}{(2m+1)^3}
	\sum_{j=\lfloor \frac{m}{2} \rfloor +5}^{3\lfloor \frac{m}{4} \rfloor} 
	\sum_{l=\lfloor \frac{m}{2} \rfloor +1}^{j-4} 
	\frac{ 1}{\sin^2 \left( \frac{\theta_{j} + \gamma_l - \pi/8}{2}\right)} 
	\sum_{\nu=1}^{\lfloor \frac{j-l}{2} \rfloor -1}
	\frac{1}{\sin \sigma_\nu }\\
	& \leq  c \sum_{j=\lfloor \frac{m}{2} \rfloor +5}^{3\lfloor \frac{m}{4} \rfloor} 
	\frac{1}{j^2}
	\sum_{l=\lfloor \frac{m}{2} \rfloor +1}^{j-4} 
	\sum_{\nu=1}^{\lfloor \frac{j-l}{2} \rfloor -1}
	\frac{1}{2\nu-1/2}\\
	& \leq c m \log(m+1).
\end{align*}
For $J(+,-)$, using the inequality $\theta_{j} + \sigma_\nu - \gamma_l > \theta_{j} - \gamma _l$,  
\begin{align*}
	J_1(+,-) & \leq \frac{1}{(2m+1)^3}
	\sum_{j=\lfloor \frac{m}{2} \rfloor +5}^{3\lfloor \frac{m}{4} \rfloor} 
	\sum_{l=\lfloor \frac{m}{2} \rfloor +1}^{j-4} 
	\frac{ 1}{\sin^2 \left( \frac{\theta_{j} -\gamma_l}{2}\right)} 
	\sum_{\nu=1}^{\lfloor \frac{j-l}{2} \rfloor -1}
	\frac{1}{\sin(\sigma_\nu)}\\
	& \leq c \sum_{j=\lfloor \frac{m}{2} \rfloor +5}^{3 \lfloor \frac{m}{4} \rfloor} 
	\sum_{l=\lfloor \frac{m}{2} \rfloor +1}^{j-4} 
	\frac{1}{\left(j-l-2\right)^2} 
	\sum_{\nu=1}^{\lfloor \frac{j-l}{2} \rfloor -1}
	\frac{1}{2\nu-1/2}\\
	& \leq c m \log(m+1).
\end{align*}
For the remaining case of $J(-,-)$, we split the sum in $\nu$,
\begin{equation*}
	\left( \sum_{\nu=1}^{ \lfloor \frac{j-l}{4} \rfloor -1} + \sum_{\nu= \lfloor \frac{j-l}{4} \rfloor}^{ \lfloor \frac{j-l}{2} \rfloor -1} \right) 
	\frac{1}{\sin(\sigma_\nu) \sin^2 \left( \frac{\theta_{j} - \sigma_\nu - \gamma_l}{2}\right)}.  
\end{equation*}
We are only considering values of $\nu \geq 1$, so we ignore any instances of $\nu=0,-1$ in the sums.  For the first sum, $\theta_{j} - \sigma_\nu - \gamma_l > (\theta_{j} - \gamma_l)/2$, so a bound of $c m\log(m+1)$ is found as in the case of $J(+,-)$.  For the second sum, 
\begin{displaymath}
	\frac{2\nu-1/2}{2m+1} \geq \frac{ (j-l-3)/2}{2m+1} \geq \frac{1}{4m+2},  	
\end{displaymath}
so it readily follows that
\begin{align*}
	J_1(-,-) & \leq c
	\sum_{j=\lfloor \frac{m}{2} \rfloor +5}^{3 \lfloor \frac{m}{4} \rfloor} 
	\sum_{l=\lfloor \frac{m}{2} \rfloor +1}^{j-4 } 
	\frac{1}{j-l-3}
	\sum_{\nu=1}^{\lfloor \frac{j-l}{2} \rfloor -1 }
	\frac{1}{\left(j-2\nu-l-1 - \tfrac{1}{4m} \right)^2} \\
	&\qquad + c m \log(m+1)\\
	&\leq c m \log(m+1).
\end{align*}
\end{proof}
The next lemma eliminates the parts of $F_{j,l}^\pm(\theta, \mp\sigma_\nu, \gamma_l)$ with $4m+2$ in the numerator.
\begin{lem} 
	\label{lowerboundlem2}
	Recalling \eqref{angledefns},
\begin{align*}
	& J_2(\pm) := \frac{1}{(2m+1)^3} 
	\sum_{j=\lfloor \frac{m}{2} \rfloor +5}^{3 \lfloor \frac{m}{4} \rfloor} 
	\sum_{l=\lfloor \frac{m}{2} \rfloor +1}^{j-4} 
	\sum_{\nu=1}^{\lfloor \frac{j-l}{2} \rfloor -1} 
	\frac{(4m+2) \sin\left(\frac{\theta_{j} \pm \sigma_\nu} {2} \right) }
	{ \sin \sigma_\nu \sin \left( \frac{ \theta_{j} \pm \sigma_\nu + \gamma_l}{2}\right) }  \leq c m\log(m+1).
\end{align*}
\end{lem}
\begin{proof}
	Since $3\pi/4 \geq \theta_{j} + \sigma_\nu + \gamma_l  \geq \theta_{j} - \sigma_\nu \geq \pi/4$, both $\sin\left(\frac{\theta_{j} \pm \sigma_\nu}{2}\right)$ and $\sin\left( \frac{\theta_{j} \pm \sigma_\nu + \gamma_l}{2} \right)$ are bounded away from zero by a positive constant.  The lemma then follows from 
\begin{align*}
	J_2(\pm) & \leq \frac{c}{2m+1}
	\sum_{j=\lfloor \frac{m}{2} \rfloor +5}^{3\lfloor \frac{m}{4} \rfloor} 
	\sum_{l=\lfloor \frac{m}{2} \rfloor +1}^{j-4} 
	\sum_{\nu=1}^{\lfloor \frac{j-l}{2} \rfloor -1} 
	\frac{1}{2\nu-1/2}\\
	&\leq c m \log(m+1).
\end{align*}
\end{proof}

By applying the triangle inequality to \eqref{lowerbound}, we obtain
\begin{align*}
	&\| \CB_{2m} \|_\infty  \geq \frac{1}{(2m+1)^3}  
	\sum_{j=\lfloor \frac{m}{2} \rfloor +1}^{3 \lfloor \frac{m}{4} \rfloor} 
	\sum_{l=\lfloor \frac{m}{2} \rfloor +1}^{j-2} 
	\sum_{\nu=1}^{\lfloor \frac{j-l}{2} \rfloor -1} \frac{1}{\sin(\sigma_\nu)}\\
	& \qquad \times \left(
	\left|F_{j,l}^+(\theta_j, \sigma_\nu, -\gamma_l) - F_{j,l}^-(\theta_j, -\sigma_\nu, -\gamma_l) \right|
	-\left| F_{j,l}^+(\theta_j, -\sigma_\nu, \gamma_l)\right| - \left| F_{j,l}^-(\theta_j, \sigma_\nu, \gamma_l) \right|
	\right).
\end{align*}
The two lemmas show that 
\begin{align*} 
	& \frac{1}{(2m+1)^3}  
	\sum_{j=\lfloor \frac{m}{2} \rfloor +1}^{3 \lfloor \frac{m}{4} \rfloor} 
	\sum_{l=\lfloor \frac{m}{2} \rfloor +1}^{j-2} 
	\sum_{\nu=1}^{\lfloor \frac{j-l}{2} \rfloor -1} \frac{1}{\sin(\sigma_\nu)}\left(
	\left| F_{j,l}^+(\theta_j, -\sigma_\nu, \gamma_l)\right| + \left| F_{j,l}^-(\theta_j, \sigma_\nu, \gamma_l) \right| \right) \\
	& \qquad \qquad \qquad \leq c m\log(m+1).
\end{align*}
We also have
\begin{align*} 
	&\left|F_{j,l}^+(\theta_j, \sigma_\nu, -\gamma_l) - F_{j,l}^-(\theta_j, -\sigma_\nu, -\gamma_l) \right|
	\geq (4m+2) 
	\left| 
	\frac{ \sin \left( \frac{\theta_{j}+\sigma_\nu}{2} \right) }
	{\sin \left( \frac{\theta_{j}-\gamma_l +\sigma_\nu}{2}\right)} -
	\frac{\sin\left( \frac{\theta_{j}-\sigma_\nu}{2}\right)}
	{\sin \left( \frac{\theta_{j} -\gamma_l -\sigma_\nu}{2}\right)}\right|\\
	&\qquad \qquad \qquad - c \left(\frac{1}{\sin^2 \left( \tfrac{\theta_j - \sigma_\nu - \gamma_l}{2} \right)} + \frac{1}{\sin^2 \left( \tfrac{\theta_j + \sigma_\nu - \gamma_l}{2} \right)} \right).
\end{align*}
Lemma \eqref{lowerboundlem1} shows that
\begin{align*} 
	&\frac{1}{(2m+1)^3}  
	\sum_{j=\lfloor \frac{m}{2} \rfloor +1}^{3 \lfloor \frac{m}{4} \rfloor} 
	\sum_{l=\lfloor \frac{m}{2} \rfloor +1}^{j-2} 
	\sum_{\nu=1}^{\lfloor \frac{j-l}{2} \rfloor -1}
	\frac{1}{\sin(\sigma_\nu)}\left(
	 \frac{1}{\sin^2 \left( \tfrac{\theta_j - \sigma_\nu - \gamma_l}{2} \right)}  +   \frac{1}{\sin^2 \left( \tfrac{\theta_j + \sigma_\nu - \gamma_l}{2} \right)} \right)\\
	& \qquad \qquad \qquad \leq c m\log(m+1).
\end{align*}
Hence, we obtain
\begin{align}
	\label{lastb2m}
	&\| \CB_{2m} \|_\infty  
	\geq \frac{1}{(2m+1)^2}  
	\sum_{j=\lfloor \frac{m}{2} \rfloor +5}^{3 \lfloor \frac{m}{4} \rfloor} 
	\sum_{l=\lfloor \frac{m}{2} \rfloor +1}^{j-2} 
	\sum_{\nu=1}^{\lfloor \frac{j-l}{2} \rfloor -1} 
	\frac{1}{\sin(\sigma_\nu)}
	\left|  
	\frac{ \sin \left( \frac{ \theta_{j}+\sigma_\nu }{ 2 } \right) } 
	{ \sin \left( \frac{\theta_{j}-\gamma_l +\sigma_\nu}{2} \right)} 
	-\frac{ \sin \left( \frac{\theta_{j}-\sigma_\nu}{2} \right) }
	{\sin \left( \frac{\theta_{j} -\gamma_l -\sigma_\nu }{2}  \right)}\right|\\
	& \qquad \qquad  - c m\log(m+1) \notag.
\end{align}

We now show the dominant part achieves the bound of  $c m (\log(m+1))^2$.  Using the formula for the product of sines, 
\begin{displaymath}
	\left| \frac{\sin\left(\frac{\theta_{j}-\sigma_\nu}{2}\right)}{\sin \left( \frac{\theta_{j}-\gamma_l -\sigma_\nu}{2}\right)}
	-\frac{\sin\left(\frac{\theta_{j}+\sigma_\nu}{2}\right)}{\sin \left( \frac{\theta_{j} -\gamma_l +\sigma_\nu}{2}\right)} \right| 
	=  \frac{\sin(\sigma_\nu) \sin\left( \frac{\gamma_l}{2}\right)}{\sin\left(\frac{\theta_{j}-\gamma_l-\sigma_\nu}{2}\right) \sin\left(\frac{\theta_{j} - \gamma_l +\sigma_\nu}{2} \right)},
\end{displaymath}
and it follows that
\begin{align*}
	& \sum_{j=\lfloor \frac{m}{2} \rfloor +1}^{3 \lfloor \frac{m}{4} \rfloor} 
	\sum_{l=\lfloor \frac{m}{2} \rfloor +1}^{j-4} 
	\sum_{\nu=1}^{\lfloor \frac{j-l}{2} \rfloor -1} 
	\frac{1}{\sin(\sigma_\nu)}
	\left| \frac{\sin\left(\frac{\theta_{j}-\sigma_\nu}{2}\right)}{\sin \left( \frac{\theta_{j} -\gamma_l -\sigma_\nu}{2}\right)}
	-\frac{\sin\left(\frac{\theta_{j}+\sigma_\nu}{2}\right)}{\sin \left( \frac{\theta_{j} -\gamma_l +\sigma_\nu}{2}\right)} \right|\\
	& \geq c \sum_{j=\lfloor \frac{m}{2} \rfloor +1}^{3 \lfloor \frac{m}{4} \rfloor} 
	\sum_{l=\lfloor \frac{m}{2} \rfloor +1}^{j-4} 
	\sum_{\nu=1}^{\lfloor \frac{j-l}{2} \rfloor -1}
	\frac{1}{(\theta_{j} - \gamma_l -\sigma_\nu)(\theta_{j} - \gamma_l +\sigma_\nu)}\\
	& = c \sum_{j=\lfloor \frac{m}{2} \rfloor +1}^{3 \lfloor \frac{m}{4} \rfloor} 
	\sum_{l=\lfloor \frac{m}{2} \rfloor +1}^{j-4} 
	\frac{1}{\theta_{j} - \gamma_l}
	\sum_{\nu=1}^{\lfloor \frac{j-l}{2} \rfloor -1}
	\left( \frac{1}{\theta_{j} - \gamma_l -\sigma_\nu} + \frac{1}{\theta_{j} - \gamma_l +\sigma_\nu} \right)\\
	& \geq  c \sum_{j=\lfloor \frac{m}{2} \rfloor +1}^{3 \lfloor \frac{m}{4} \rfloor} 
	\sum_{l=\lfloor \frac{m}{2} \rfloor +1}^{j-4}
	\frac{1}{\theta_{j} - \gamma_l}
	\sum_{\nu=1}^{\lfloor \frac{j-l}{2} \rfloor -1}
	\frac{1}{\theta_{j} - \gamma_l -\sigma_\nu}.
\end{align*}
Dividing by $(2m+1)^2$, we have
\begin{align*}
	&\frac{1}{(2m+1)^2} 
	\sum_{j=\lfloor \frac{m}{2} \rfloor +1}^{3 \lfloor \frac{m}{4} \rfloor} 
	\sum_{l=\lfloor \frac{m}{2} \rfloor +1}^{j-4} 
	\frac{1}{\theta_{j}-\gamma_l}
	\sum_{\nu=1}^{\lfloor \frac{j-l}{2} \rfloor -1}
	\frac{1}{\theta_{j}-\sigma_\nu-\gamma_l}\\
	& \geq c 
	\sum_{j=\lfloor \frac{m}{2} \rfloor +1}^{3 \lfloor \frac{m}{4} \rfloor} 
	\sum_{l=\lfloor \frac{m}{2} \rfloor +1}^{j-4} 
	\frac{1}{j-(2l+1)(\frac{2m+1}{4m})}
	\sum_{\nu=1}^{\lfloor \frac{j-l}{2} \rfloor -1}
	\frac{1}{j-2\nu-(2l+1)(\frac{2m+1}{4m})}\\
	&\geq c \sum_{j=\lfloor \frac{m}{2} \rfloor +1}^{3 \lfloor \frac{m}{4} \rfloor} 
	\sum_{l=\lfloor \frac{m}{2} \rfloor +1}^{j-4} 
	\frac{1}{j-l-1/2}
	\sum_{\nu=1}^{\lfloor \frac{j-l}{2} \rfloor -1 }
	\frac{1}{j-2\nu -l-1/2}\\
	&\geq c \sum_{j=\lfloor \frac{m}{2} \rfloor +1}^{3 \lfloor \frac{m}{4} \rfloor} 
	\sum_{l=\lfloor \frac{m}{2} \rfloor +1}^{j-4} 
	\frac{\log\left( j-l \right)}{j-l-1/2}\\
	& \geq c m (\log(m+1))^2,
\end{align*}
which completes the proof of the lower bound.

\subsection{Upper Bound}
As mentioned previously, the estimate of the Lebesgue constant on $B^2 \times [-1,1]$ does not trivially reduce to estimates on $B^2$ and $[-1,1]$.  Moreover, a straight-forward estimate of the upper bound of the Lebesgue constant, as done in \cite{XuRadon} for the case on $B^2$, would be extremely lengthy and cumbersome.  We instead use a different approach, by deriving generating functions for the series in the definition of $T_{\nu,j,l}(x,y,z)$, and then writing $T_{\nu,j,l}(x,y,z)$ as the Fourier coefficient of the product of the generating functions.  The idea for this approach comes from \cite{LiXu}, and provides an alternative proof for the upper bound of \cite[Theorem 5.2]{XuRadon}. 

\begin{lem} 
	\label{lemmaobscene}
	For $0 < |r|< 1$ and $m \geq 0$, 
	\begin{multline*}
		T_{\nu,j,l}(x,y,z) = \frac{1}{(2m+1)^3} \frac{1}{\sin \sigma_\nu(x,y) } \\
		\frac{1}{2\pi} 
		\int_0^{2\pi} \, \frac{1}{1-re^{i\theta}} G_1(re^{i\theta}, z, z_l) G_2(re^{i\theta}, \theta_j, \sigma_\nu(x,y))e^{-2mi\theta}  d\theta r^{-2m},
	\end{multline*}
	where
	\begin{displaymath}
		G_1(r,z,z_l):=\sum_{k=0}^\infty \tilde{T}_k(z) \tilde{T}_k(z_l) r^k,
	\end{displaymath}
	and
	\begin{displaymath}
		G_2(r,\theta_j,\sigma_\nu(x,y)):= \sum_{k=0}^\infty (k+1) \sin((k+1)(\theta_{j})) \sin((k+1)(\sigma_\nu(x,y))) r^k.
	\end{displaymath}
\end{lem}
\begin{proof}
Our first step is to derive the generating function of the function
	\begin{multline*}
		R_{N}(\theta_j, \sigma_\nu(x,y),z,z_l)\\
		:=\sum_{n=0}^{N} \sum_{k=0}^n (k+1) \sin( (k+1)\theta_j ) \sin( ( k+1) \sigma_\nu(x,y) ) 
		\tilde{T}_{n-k}(z_l) \tilde{T}_{n-k}(z)
	\end{multline*}
in $T_{\nu,j,l}(x,y,z)$.  Since the coefficient of $r^{N}$ in
\begin{displaymath}
	\sum_{n=0}^\infty \sum_{k=0}^n \sin( (k+1) \theta_j ) \sin( (k+1) \sigma_\nu(x,y) ) r^n 
	\sum_{j=0}^\infty \tilde{T}_j(z) \tilde{T}_j(z_l) r^j
\end{displaymath} 
is precisely $R_{N}(\theta_j,\sigma_\nu(x,y),z,z_l)$, and
\begin{align*}
	\frac{1}{1-r} \sum_{k=0}^\infty \sin( (k+1) \theta_j ) &\sin( (k+1) \sigma_\nu(x,y) ) r^k \\
	& = \sum_{n=0}^\infty \sum_{k=0}^n \sin( (k+1) \theta_j ) \sin( (k+1) \sigma_\nu(x,y) ) r^n ,
\end{align*}
it follows that 
\begin{displaymath}
	\sum_{N=0}^\infty R_N(\theta_j, \sigma_\nu(x,y)re^, z, z_l) r^k = \frac{1}{1-r} G_1(r,z,z_l)G_2(r,\theta_j,\sigma_\nu(x,y)).
\end{displaymath}
Since both sides of the above equation are analytic functions of $r$ for $|r|<1$, we may replace $r$ with $re^{i\theta}$ to obtain analytic, complex-valued functions of $r$.  Since
\begin{displaymath}
	R_{2m}( \theta_j, \sigma_\nu(x,y),z,z_l) = \frac{1}{2\pi} \int_0^{2\pi} \, 
	\sum_{N=0}^\infty R_N(\theta_j, \sigma_\nu(x,y), z, z_l) r^{2m}e^{2mi\theta} e^{-2mi\theta} d\theta r^{-2m},
\end{displaymath}
it follows that
\begin{multline*}
	R_{2m}(\theta_j, \sigma_\nu(x,y), z, z_l) \\
	= \frac{1}{2\pi} \int_0^{2\pi} \,
	\frac{1}{1-re^{i\theta}} G_1(re^{i\theta}, z, z_l) G_2(re^{i\theta}, \theta_j, \sigma_\nu(x,y)) e^{-2mi\theta} 
	d\theta \, r^{-2m}.
\end{multline*}
The lemma follows from the fact that 
\begin{displaymath}
	(2m+1)^3 \sin \sigma_\nu(x,y) R_{2m}(\theta_j, \sigma_\nu(x,y), z, z_l) = T_{\nu,j,l}(x,y,z).
\end{displaymath}
\end{proof}

We next obtain compact formulas for $G_1(re^{i\theta}, z, z_l)$ and $G_2(re^{i\theta}, \theta_j, \sigma_\nu(x,y))$, and obtain estimates for these functions.

\begin{lem}  For $m \geq 0$, and $r=1-\frac{1}{m}$,
\begin{align}
	\label{B2Mestimate}
	\left|T_{\nu,j,l}(x,y,z)\right| 
	&  \leq \frac{1}{(2m+1)^3} \frac{1}{\sin \sigma_\nu(x,y)} \frac{1}{2\pi}
	\int_0^{2\pi} \frac{1}{\left| 1- re^{i\theta}\right|} \\
	& \notag \quad \times \left( \left| A_1^+(x,y) - A_1^-(x,y) \right| + \left| A_2^+(x,y) - A_2^-(x,y) \right| \right)\\
	& \notag \quad \times \left( \left|P(re^{i\theta}, \gamma_z + \gamma_l)\right| + \left|P(re^{i\theta}, \gamma_z - \gamma_l)\right| \right) d\theta.
\end{align}
where
\begin{equation}
	A_1^{\pm}(x,y)
	=\frac{1}{1-2re^{i\theta}\cos\left( \theta_{j} \pm \sigma_\nu(x,y) \right) + r^2e^{2i\theta}},
	\label{A1}
\end{equation}
and
\begin{equation}
	A_2^{\pm}(x,y) 
	=\frac{(1-r^2e^{2i\theta})(re^{i\theta} - \cos\left( \theta_{j} \pm \sigma_\nu(x,y) \right))}
	{(1-2 r e^{i\theta} \cos \left( \theta_{j} \pm \sigma_\nu(x,y) \right) + r^2 e^{2 i \theta} )^2}.
	\label{A2}
\end{equation}
\end{lem}
\begin{proof}
	First, it follows from Lemma \eqref{lemmaobscene} that 
	\begin{multline*}
		\left| T_{\nu,j,l} (x,y,z) \right| \leq \frac{1}{(2m+1)^3} \frac{1}{\sin \sigma_\nu(x,y)} \\
		\times \frac{1}{2\pi} \int_0^{2\pi} \frac{1}{|1-re^{i\theta}|} \left| G_1(re^{i\theta},z,z_l) \right|
		\left| G_2(re^{i\theta}, \theta_j, \sigma_\nu(x,y)) \right| d\theta r^{-2m}.
	\end{multline*}
	We next derive estimates for $|G_1(re^{i\theta}, z, z_l)|$ and $|G_2(re^{i\theta}\theta_j, \sigma_\nu(x,y)|$.  The compact formula for the generating function $G_1(re^{i\theta},z,z_l)$ is well-known,
\begin{equation}
	\label{G1compact}
	G_1(re^{i\theta},z,z_l)=\frac{1}{4} \left[ P(re^{i\theta},\gamma_z + \gamma_l) + P(re^{i\theta},\gamma_z - \gamma_l) \right], 
\end{equation}
where $P(r,\phi)$ is the Poisson kernel, defined by
	\begin{equation}
		P(r,\phi)
		:= 1+ 2\sum_{n=1}^\infty \cos(n\phi) r^n 
		= \frac{1-r^2}{1 - 2 r \cos \phi + r^2},
		\ \ 0 \leq r \leq 1.
		\label{Poisson}
	\end{equation}
For $G_2(re^{i\theta},\theta_j,\sigma_\nu(x,y))$, we use the identity for the product of sines to obtain
\begin{align*}
	&G_2(re^{i\theta},\theta_j,\sigma_\nu)\\ 
	& \qquad = \frac{1}{2} \frac{d}{dr} \sum_{k=1}^\infty \left[ \cos( (k) ( \theta_j - \sigma_\nu(x,y) ) ) -  \cos((k+1)(\theta_{j}+\sigma_\nu(x,y) ) \right]  r^ke^{ik\theta}\\
	 & \qquad =\frac{1}{8} \frac{d}{dr}  \left[P(re^{i\theta}, \theta_j - \sigma_\nu(x,y)) - P(re^{i\theta},\theta_j + \sigma_\nu(x,y))\right].
\end{align*}
Using the formula
\begin{displaymath}
	\frac{d}{dr} P(r,\phi) 
	= -2 \left( \frac{r}{1-2r\cos\phi + r^2} + \frac{(r-\cos \phi)(1-r^2)}{(1-2r\cos\phi + r^2)^2} \right),
\end{displaymath}
we obtain
\begin{align*}
	G_2&(re^{i\theta},\theta_j,\sigma_\nu)=-\frac{1}{4}\\
	\times&\bigg[ \frac{re^{i\theta}}{1-2re^{i\theta}\cos(\theta_j - \sigma_\nu(x,y)) + r^2e^{2i\theta}}
	+\frac{[re^{i\theta}-\cos(\theta_j - \sigma_\nu(x,y))](1-r^2e^{2i\theta})}{ \left[ 1 - 2re^{i\theta} \cos( \theta_j - \sigma_\nu(x,y)) + r^2e^{2i\theta} \right]^2}\\
	-&\frac{re^{i\theta}}{1-2re^{i\theta}\cos(\theta_j + \sigma_\nu(x,y)) + r^2e^{2i\theta}} 
	+\frac{[re^{i\theta}-\cos(\theta_j + \sigma_\nu(x,y))](1-r^2e^{2i\theta})}{ \left[ 1 - 2re^{i\theta}\cos( \theta_j + \sigma_\nu(x,y)) + r^2e^{2i\theta} \right]^2}\bigg].
\end{align*}

It follows that 
\begin{displaymath}
	|G_2(re^{i\theta}, \theta_j, \sigma_\nu(x,y))| \leq  \left( \left| A_1^+(x,y) - A_1^-(x,y) \right| + \left| A_2^+(x,y) - A_2^-(x,y) \right| \right).
\end{displaymath}
Finally, since the inequality in \eqref{lemmaobscene} holds for all values of $r$ with $0 < r < 1$, we set $r = 1 - \frac{1}{m}$.  With this choice of $r$, $r^{-2m}$ converges to $e^{-2}$ as $m$ approaches infinity, and so $r^{-2m}$ is bounded by a constant for all $m$.    
\end{proof}
Before beginning the estimate, we make several reductions in the range of the sums and values of $x,y,z$ that need to be considered.  
\begin{enumerate}
	\item First, we can reduce the interval of integration to $[0,\pi]$. To see this, replace $\theta$ with $2\pi - \theta$.  This change of variable amounts to conjugation of the complex number $re^{i\theta}$, and hence the norms of the expression are unchanged.  

	\item We may also restrict $\gamma_z$ to the interval $[0,\pi/2m]$.  To see this, replace $\gamma_z$ with $\gamma_z + \tfrac{\pi}{2m}$ in $P(re^{i\theta}, \gamma_z +\theta_j)$ and $P(re^{i\theta}, \gamma_z - \gamma_l)$.  We see that, upon changing the summation index from $l$ to $2m-1 -l$, 
\begin{align*}
	&\sum_{l=0}^{2m-1} \frac{1-r^2e^{2i\theta}}{1-2re^{i\theta}\cos \left( \gamma_z + \frac{\pi}{2m} + \gamma_l \right) + r^{2}e^{2i\theta}} \\
	& = \sum_{l=1}^{2m-1} \frac{1-r^2e^{2i\theta}}{1-2re^{i\theta}\cos \left( \gamma_z + \gamma_l \right) + r^{2}e^{2i\theta} } + \frac{1-r^2 e^{2i\theta}}{1-2re^{i\theta}\cos\left( \gamma_z - \frac{4m-1}{4m}\pi \right) + r^2 e^{2i\theta}},
\end{align*}
and
\begin{align*}
	&\sum_{l=0}^{2m-1} \frac{1-r^2e^{2i\theta}}{1-2re^{i\theta}\cos \left( \gamma_z + \frac{\pi}{2m} - \gamma_l \right) + r^{2}e^{2i\theta}} \\
	& \qquad = \sum_{l=0}^{2m-2} \frac{1-r^2e^{2i\theta}}{1-2re^{i\theta}\cos \left( \gamma_z - \gamma_l \right) + r^{2}e^{2i\theta} } + \frac{1-r^2 e^{2i\theta}}{1-2re^{i\theta}\cos\left( \gamma_z + \frac{\pi}{4m} \right) +r^2 e^{2i\theta}}.
\end{align*}
It follows that the expression
\begin{displaymath}
	\sum_{l=0}^{2m-1} \left| P(re^{i\theta}, \gamma_z + \gamma_l) + P(re^{i\theta}, \gamma_z - \gamma_l) \right|
\end{displaymath}
is invariant under translations of $\gamma_z$ by $\pi/2m$, so we only need to consider $\gamma_z \in \left[ 0, \tfrac{\pi}{2m} \right]$.

\item The sum in $j$ may be reduced to $1\leq j \leq m$.  Replacing $j$ with $2m+1-j$, $\sin \theta_{2m+1-j} =\sin \theta_j $, and $\cos(\theta_{2m+1-j} \pm \sigma_\nu(x,y) ) = \cos(\theta_j \mp (\pi-\sigma_\nu(x,y)) )$.  It follows from the definition of $\sigma_\nu(x,y)$ that $\pi-\sigma_\nu(x,y) = \sigma_\nu(-x,-y)$, which implies that $\sin \sigma_\nu(-x,-y)  = \sin \sigma_\nu(x,y)$.  Hence, 
\begin{align*}
	&\sum_{j=m+1}^{2m} \frac{\sin \theta_j}{\sin \sigma_\nu(x,y)} \left(\left|A_1^+(x,y) - A_1^- (x,y) \right| + \left| A_2^+(x,y) - A_2^-(x,y)\right|\right)\\
	&=\sum_{j=1}^m \frac{\sin \theta_j}{\sin \sigma_\nu(-x,-y)} \\
	&\qquad \qquad  \times \left(\left|A_1^+(-x,-y) - A_1^- (-x,-y) \right| + \left| A_2^+(-x,-y) - A_2^-(-x,-y)\right|\right),
\end{align*}
which shows we only need to consider $1 \leq j \leq m$.

\item We also only need to consider $(x,y)$ in the region 
\begin{displaymath}
	\Gamma_m=\left\{ (\rho \cos \eta,  \rho \sin \eta )\ :\ -\frac{\pi}{4m+2} \leq \eta \leq \frac{\pi}{4m+2} \right\}.
\end{displaymath}
To see this, let $x=\rho \cos \eta $ and $y=\rho \sin \theta $, so that $\cos \sigma_\nu(x,y) =\rho \cos(\eta-\phi_\nu)$.  Note that the collection of points $\rho \cos(\eta - \phi_\nu)$, for $\nu=0,1,\ldots,2m$, is unchanged by a rotation of $\eta$ by $\phi_\nu$.  Moreover, every expression involving $\sigma_\nu(x,y)$ in the right side of \eqref{B2Mestimate} can be written in terms of $\cos \sigma_\nu(x,y)$.  Since $ 0 \leq \sigma_\nu(x,y) \leq \pi$, $\sin \sigma_\nu(x,y) = \sqrt{1 - \cos^2 \sigma_\nu(x,y)}$, and the expressions $\cos(\theta_j \pm \sigma_\nu(x,y))$ can be expanded using the cosine addition identity.  Hence, every expression involving $\sigma_\nu(x,y)$ in 
\begin{displaymath}
	\sum_{\nu=0}^{2m} \frac{1}{\sin \sigma_\nu (x,y)} \left( \left| A_1^+(x,y) - A_1^-(x,y) \right| + \left|A_2^+(x,y)- A_2^-(x,y) \right| \right)
\end{displaymath}
is the same at the points $(\rho \cos \eta, \rho \sin \eta)$ and $(\rho \cos (\eta+\phi_\nu), \rho \sin (\eta + \phi_\nu))$.

\item Finally, we may also reduce the sum in $\nu$ to $ 0 \leq \nu \leq m$.  First note that
\begin{displaymath}
	\cos \sigma_{2m+1-\nu}(x,y)  = \cos \sigma_\nu(x,-y) , \text{ and } \sin \sigma_{2m+1-\nu}(x,y)  = \sin \sigma_\nu(x,-y).
\end{displaymath}
Hence, we obtain 
\begin{align*}
	&\sum_{\nu=m+1}^{2m} \frac{1}{\sin \sigma_\nu (x,y)} \left( \left| A_1^+(x,y) - A_1^-(x,y) \right| + \left|A_2^+(x,y)- A_2^-(x,y) \right| \right)\\
	=&\sum_{\nu=1}^{m} \frac{1}{\sin \sigma_\nu (x,-y)} \left( \left| A_1^+(x,-y) - A_1^-(x,-y) \right| + \left|A_2^+(x,-y)- A_2^-(x,-y) \right| \right),
\end{align*}
and since $\Gamma_m$ is symmetric with respect to $y$, we only need to consider $0 \leq \nu \leq m$.
\end{enumerate}

From these reductions, it follows that
\begin{align*}
	&\|\mathcal{B}_{2m}\|_\infty \leq c \max_{ \substack{
			(x,y) \in \Gamma_m\\
			\gamma_z \in \left[ 0, \frac{\pi}{2m} \right] }}  
		\int_0^{2\pi} 
		\frac{1}{2m+1} \sum_{l=0}^{2m-1} 
		\left|P(re^{i\theta},\gamma_z + \gamma_l)  +  P(re^{i\theta},\gamma_z-\gamma_l)  \right|\\
		& \notag \qquad \times 
		\frac{1}{(2m+1)^2}
		\sum_{\nu=0}^{2m} \sum_{j=1}^{2m} 
		\frac{1}{\left| 1- re^{i\theta} \right|}
		\frac{ \sin \theta_{j} }{ \sin \sigma_\nu(x,y)  }\\ 
		&\notag \qquad \times \left( \left| A_1^+(x,y) - A_1^-(x,y) \right| + \left| A_2^+(x,y) - A_2^-(x,y) \right|\right) 
		d\theta .
\end{align*}

Consequently, the proof of Theorem \ref{Lebesgueconstant} will follow from the three lemmas below. 
\begin{proof}[]
\end{proof}
\begin{lem}
	\label{lem:upperboundz}
	For $z \in [0,\pi/2m]$,
\begin{displaymath}
	\frac{1}{2m+1} \sum_{l=0}^{2m-1} \left|P(re^{i\theta},\gamma_z + \gamma_l)\right|  + \left| P(re^{i\theta},\gamma_z-\gamma_l)  \right| \leq c \log(m+1)
\end{displaymath}
for some $c$ which is independent of $\theta$, $z$, and $m$.  
\end{lem}

\begin{lem}
	\label{lem:upperboundA1}
	For $(x,y) \in \Gamma_m$, 
\begin{align*}
	&\frac{1}{(2m+1)^2} \int_{0}^\pi \frac{1}{|1-re^{i\theta}|} \sum_{\nu=0}^{2m} \sum_{j=1}^{2m} 
	\frac{\sin(\theta_{j,2m})}{\sin(\sigma_\nu(x,y))} \left| A_1^+(x,y) - A_1^-(x,y) \right| d\theta \leq c m \log(m+1)
\end{align*}
for some $c$ which is independent of $x,y$, and $m$.
\end{lem}

\begin{lem}
	\label{lem:upperboundA2}
	For $(x,y) \in \Gamma_m$,
\begin{align*}
	&\frac{1}{(2m+1)^2} 
	\int_{0}^\pi 
	\frac{1}{|1-re^{i\theta}|} 
	\sum_{\nu=0}^{2m} \sum_{j=1}^{2m} 
	\frac{ \sin(\theta_{j,2m}) }{\sin(\sigma_\nu(x,y))} 
	\left| A_2^+(x,y) - A_2^-(x,y) \right|
	d\theta \leq c m (\log(m+1))
\end{align*}
for some $c$ which is independent of $x,y$, and $m$.
\end{lem}

The proofs of these lemmas are contained in the next section.

\section{Proofs of Lemmas \ref{lem:upperboundz}, \ref{lem:upperboundA1}, and \ref{lem:upperboundA2}}
\setcounter{equation}{0}

\subsection{Proof of Lemma \ref{lem:upperboundz}}

\begin{proof}
Fix $z \in \left[0,\frac{\pi}{2m}\right]$.  First applying the triangle inequality, we have
\begin{align*}
	&\left| G_1(r,\theta,z)\right|\\
	&\qquad \leq c \left( 
	\left| \frac{ 1 -r^2e^{2i\theta} }{1-2re^{i\theta} \cos ( \gamma_z + \gamma_l) + r^2e^{2i\theta}}\right|
	+ \left| \frac{ 1 -r^2e^{2i\theta} }{1-2re^{i\theta} \cos ( \gamma_z - \gamma_l) + r^2e^{2i\theta}} \right|
	\right).
\end{align*}
Next we apply the approximation
\begin{equation}
	\left| 1-re^{i\theta} \right| \asymp \left| \sin \tfrac{\theta}{2} \right| + m^{-1}
	\label{easyapprox}
\end{equation}
and also note that
\begin{equation}
	1-2re^{i\theta}\cos(\phi) + r^2e^{2i\theta]} = (1-re^{i(\theta+\phi)})(1-re^{i(\theta-\phi)})
	\label{easyfactor}
\end{equation}
to obtain 
\begin{align*}
	\left| G_1(r,\theta,z)\right|&
	\leq c\left( 
	\frac{ 
	\left( \left| \sin\left( \frac{\pi-\theta}{2}\right)\right| + m^{-1} \right)
	\left( \left| \sin \left( \frac{\theta}{2} \right)\right| + m^{-1} \right)}
	{\left( \left|\sin \left( \frac{\theta+\gamma_z+\gamma_l}{2} \right)\right| + m^{-1} \right)
	\left( \left| \sin \left( \frac{\theta -\gamma_z -\gamma_l}{2}\right)\right| + m^{-1} \right)}\right.\\
	&\qquad 
	\left. + 
	\frac{ 
	\left( \left| \sin\left( \frac{\pi-\theta}{2} \right) \right| + m^{-1} \right)
	\left( \left| \sin \left( \frac{\theta}{2} \right) \right| + m^{-1} \right) }
	{\left( \left|\sin \left( \frac{\theta+\gamma_z-\gamma_l}{2} \right)\right| + m^{-1} \right)
	\left( \left| \sin \left( \frac{\theta -\gamma_z +\gamma_l}{2}\right)\right| + m^{-1} \right)}
	\right)\\
	& := c(\Psi_1 + \Psi_2)
\end{align*}
\textbf{Case 1: $0 \leq \theta \leq \pi/2$.}

We may ignore the factor of $\sin \left( \frac{\pi-\theta}{2}\right)$ in this case.  We also have
\begin{displaymath}
	\left| \frac{ \theta+\theta_z + \gamma_l}{2}\right| \leq \frac{3\pi}{4} + \frac{\pi}{4m},
\end{displaymath}
so that $\sin( (\theta\pm_1 \gamma_z \pm_2 \gamma_l)/2 ) \approx \theta \pm_1 \gamma_z \pm_2 \gamma_l$.  Fixing $\theta$, for $\Psi_1$, we have 
\begin{align*}
	\frac{1}{2m+1} \sum_{l=0}^{2m-1} \Psi_1 
	&\leq \frac{c}{2m+1} \sum_{l=0}^{2m-1} \frac{\theta + m^{-1}}{ (\theta+\gamma_z+\gamma_l + m^{-1})(\left| \theta - \gamma_z - \gamma_l\right| + m^{-1})}\\
	& \leq \frac{c}{2m+1} \sum_{l=0}^{2m-1} \frac{1}{\left| \theta - \gamma_z - \gamma_l \right| + m^{-1}}\\
	& \leq c \sum_{l=1}^{2m} \frac{1}{\left| m(\theta - \gamma_z) - \frac{2l+1}{4} \pi \right| + 1}\\ 
	& \leq c \sum_{l=0}^{2m} \frac{1}{l+1} \leq c \log(m+1).
\end{align*}
This type of estimate will be very common throughout the proof, and we will frequently omit the repetitive details.

For $\Psi_2$, we first sum starting from $l=1$, so that $\theta + \gamma_l - \gamma_z > \theta$, and hence
\begin{align*}
	&\frac{1}{2m+1} \sum_{l=1}^{2m-1} \frac{\left|\theta\right| +m^{-1} }
	{\left(\left| \theta + \gamma_z - \gamma_l\right| + m^{-1}\right)(\left( \left| \theta-\gamma_z + \gamma_l \right| + m^{-1} \right)}\\
	\leq & \frac{1}{2m+1} \sum_{l=1}^{2m-1} \frac{1}{ \left| \theta + \gamma_z - \gamma_l \right| + m^{-1}}
	\leq  c \log(m+1).
\end{align*}
For the term corresponding to $l=0$, either $\gamma_z \geq \pi/4m$ or $\gamma_z \leq \pi/4m$; assuming, without loss of generality, the former, it then follows that
\begin{align*}
	&\frac{1}{2m+1} \frac{\theta+m^{-1}}{\left(\left|\theta+\gamma_z-\frac{\pi}{4m}\right| + m^{-1}\right) 
	\left( \left| \theta-\gamma_z + \frac{\pi}{4m}\right| + m^{-1} \right)}\\
	& \leq \frac{1}{2m+1} \frac{1}{ \left| \theta - \gamma_z + \frac{\pi}{4m} \right| + m^{-1} }
	\leq 1,
\end{align*}
which shows that 
\begin{displaymath}
	\frac{1}{2m+1} \sum_{l=0}^{2m-1} \Psi_2 \leq c \log(m+1) 
\end{displaymath}
for the case when $\theta \in [0,\pi/2]$.

\textbf{Case 2: $\pi/2 \leq \theta \leq \pi$.}
For this case, the factor of $\sin(\theta/2)$ can be ignored.  For $\Psi_1$, note that
\begin{align*}
	&-\frac{\pi}{4m} \leq \frac{2\pi-\theta-\gamma_z-\gamma_l}{2} \leq \frac{3\pi}{4},\ 
	 -\frac{\pi}{4} - \frac{\pi}{8m} \leq \frac{\theta - \gamma_z -\gamma_l}{2} \leq \frac{\pi}{2},
\end{align*}
so we may approximate the sine functions accordingly and obtain
\begin{align*}
	\frac{1}{2m+1} \sum_{l=0}^{2m-1} \Psi_2 & 
	\leq \frac{c}{2m+1} \sum_{l=0}^{2m-1} \frac{ \left| \pi - \theta \right| + m^{-1}}
	{\left( \left| 2\pi - \theta - \gamma_z -\gamma_l\right| + m^{-1}\right) \left( \left| \theta - \gamma_z -\gamma_l \right| + m^{-1} \right) }\\ 
	& = \frac{c}{2m+1} \sum_{l=0}^{2m-1} \frac{\left|\pi-\theta\right|+m^{-1}}
	{\left( \left| \pi-\theta-\gamma_z+\gamma_l \right| +m^{-1}\right)
	 \left( \left| \pi-\theta + \gamma_z -\gamma_l \right| + m^{-1} \right) },
\end{align*}
where we substituted $2m-1-l$ for $l$ in the last equality.  This estimate is very similar to the estimate of $\Psi_2$ in Case 1, and hence for $\pi/2\leq \theta \leq  \pi$,
\begin{displaymath}
	\frac{1}{2m+1} \sum_{l=1}^{2m} \Psi_1\leq c \log (m+1).
\end{displaymath}

For $\Psi_2$,
\begin{align*}
	& 0 \leq \frac{2\pi-\theta+\gamma_z-\gamma_l}{2}\leq \frac{3\pi}{4} + \frac{\pi}{4m},\ 
	 -\frac{\pi}{2} \leq \frac{\theta+\gamma_z-\gamma_l}{2} \leq \frac{\pi}{2} + \frac{\pi}{4m},
\end{align*}
and so 
\begin{align*}
	\frac{1}{2m+1} \sum_{l=0}^{2m-1} \Psi_2 & 
	\leq \frac{c}{2m+1} \sum_{l=0}^{2m-1} \frac{ \left| \pi-\theta \right|+m^{-1} }
	{ \left( \left| 2\pi - \theta + \gamma_z - \gamma_l \right| + m^{-1} \right)
	  \left( \left| \theta + \gamma_z -\gamma_l \right| + m^{-1} \right) }\\
	  & = \frac{c}{2m+1} \sum_{l=0}^{2m-1} \frac{ \left| \pi - \theta \right| + m^{-1} }
	{ \left( \left| \pi-\theta+\gamma_z + \gamma_l \right| + m^{-1} \right)
	\left( \left| \pi - \theta - \gamma_z - \gamma_l \right| + m^{-1}\right) },
\end{align*}
which is very similar to our estimate for $\Psi_1$ in Case 1, and we again get a estimate of $c \log(m+1)$.
\end{proof}

\subsection{Proof of Lemma \ref{lem:upperboundA1}}

\begin{proof}
	The proof of the remaining two lemmas will proceed by separating the integral into three different regions, then dividing the sums in $j$ and $\nu$ into several sections, and performing estimates on each resulting section.  Frequently, obtaining an estimate consists of bounding quotients by a constant, and then estimating similar ``types'' of sums and integrals.  For the sake of brevity, we list these types here, and then direct the reader to the type of estimate that arises in the each piece.  The symbols $\phi_1, \phi_2$ and $\xi$ refer to values that are particular to the section under investigation.  For two different expressions $f_1$ and $f_2$, the notation $\{ f_1, f_2\}$ indicates that either expression satisfies that type.  Finally, we note that these estimates also hold for sums whose range of indices are a subset of those listed below.   

\noindent
\textbf{Type 1}:  
\begin{align*}
	&\frac{1}{(2m+1)^2} \int_{\phi_1}^{\phi_2} \frac{1}{ \left\{ \left| \theta-\phi_1 \right|,\,  \left| \phi_2 - \theta \right| \right\} + m^{-1}}
	\sum_{\nu=0}^m \sum_{j=1}^m \frac{1}{\left| \theta_{j} + \xi \right| + m^{-1} } \ d\theta \\
	& \qquad \qquad \leq c(\log(m+1))^2.
\end{align*}

\noindent
\textbf{Type 2}:
\begin{align*}
	&\frac{1}{(2m+1)^2} \int_{\phi_1}^{\phi_2} \sum_{\nu=0}^m \sum_{j=1}^m \frac{1}{( \left| \theta_{j} + \xi \right| + m^{-1})^2} \ d\theta \leq cm.
\end{align*}

\noindent
\textbf{Type 3}:
\begin{align*}
	& \frac{1}{(2m+1)^2} \int_{\phi_1}^{\phi_2} 
	\frac{1}{\left( \left\{ \left| \theta-\phi_1 \right|, \, \left| \phi_2 - \theta \right| \right\} + m^{-1} \right)^2}
	\sum_{\nu=0}^m \sum_{j=1}^m \frac{1}{\left| \theta_{j} + \xi \right| + m^{-1} } \ d\theta \\
	& \qquad \qquad \leq c m \log(m+1).
\end{align*}

\noindent
\textbf{Type 4}:
\begin{align*}
	& \frac{1}{(2m+1)^2} \int_{\phi_1}^{\phi_2} 
	\frac{1}{ \left\{ \left| \phi_2 - \theta \right|, \left| \theta-\phi_1 \right| \right\} + m^{-1} }
	\sum_{\nu=0}^m  \sum_{j=1}^m  \frac{1}{\left( \left| \theta_{j} + \xi \right| + m^{-1} \right)^2} \ d\theta \\
	& \qquad \qquad \leq c m \log(m+1).
\end{align*}

\noindent
\textbf{Type 5}:
\begin{align*}
	& \frac{1}{(2m+1)^2} \int_{\phi_1}^{\phi_2} 
	\frac{1}{ \left\{ \left| \phi_2 - \theta \right|, \left| \theta-\phi_1 \right| \right\} + m^{-1} }
	\sum_{\nu=0}^m  \frac{1}{\left\{ \left| \sigma_\nu(x,y) \right|, \, \left| \pi-\sigma_\nu(x,y) \right| \right\} + m^{-1} }\\
	& \qquad \qquad \times \sum_{j=1}^m  \frac{1}{\left| \theta_{j} + \xi \right| + m^{-1}  } \ d\theta \, 
	\leq c (\log(m+1))^3.
\end{align*}

\noindent
\textbf{Type 6}:
\begin{align*}
	& \frac{1}{(2m+1)^2} \int_{\phi_1}^{\phi_2} 
	\sum_{\nu=0}^m \frac{1}{ \left\{ \left| \sigma_\nu(x,y) \right|,\, \left| \pi-\sigma_\nu(x,y) \right| \right\} + m^{-1} }
	\sum_{j=1}^m  \frac{1}{\left| \theta_{j} + \xi \right| + m^{-1} }\ d\theta\\
	& \qquad \qquad \leq c  (\log(m+1))^2.
\end{align*}

\noindent
\textbf{Type 7}:
\begin{align*}
	& \frac{1}{(2m+1)^2} \int_{\phi_1}^{\phi_2}
	\sum_{\nu=0}^m \frac{1}{\left\{ \left| \sigma_\nu(x,y) \right|, \, \left| \pi-\sigma_\nu(x,y) \right| \right\} + m^{-1} }
	\sum_{j=1}^m \frac{1}{ \left( \left| \theta_{j} + \xi \right| + m^{-1} \right)^2} \ d\theta\\
	& \qquad \qquad \leq c m \log(m+1) .
\end{align*}

The above estimates are easily obtained with the following proposition and lemma.

\begin{prop} For any real number $\xi$, 
\begin{displaymath}
	\sum_{j=1}^m \frac{1}{(\left| \theta_{j} + \xi \right| + m^{-1})^k} \leq 
	\left\{ \begin{array}{ll}
		c (2m+1) \log(m+1) \ & \ k=1\\
		c  (2m+1)^2 \ & \ k=2\\
	\end{array}
	\right.
\end{displaymath}
\end{prop}

\begin{proof}
By the definition of $\theta_{j}$,
\begin{align*}
	\sum_{j=1}^m \frac{1}{\left( \left| \theta_{j} + \xi \right| + m^{-1} \right)^k}
	&=  (2m+1)^k \sum_{j=1}^{m} \frac{1}{\left| j\pi + (2m+1)\xi \right| + 2 + \frac{1}{m}}\\
	&\leq (2m+1)^k \, 2 \,\sum_{j=1}^m \frac{1}{(j + 1)^k},
\end{align*}
from which the proposition easily follows.
\end{proof}

\begin{lem}
For $(x,y) \in \Gamma_m$,
\begin{align}
	& \sum_{\nu=0}^m \frac{1}{\left| \sigma_\nu(x,y) \right| + m^{-1}} \leq c m \log(m+1)
	\label{sigmaest}\\
	& \sum_{\nu=0}^m \frac{1}{\left| \pi - \sigma_\nu(x,y) \right| + m^{-1}  } \leq c m \log(m+1).
	\label{pisigmaest}
\end{align}
\end{lem}
\begin{proof}
Recall that $(x,y) =  (r\cos(\phi), r\sin(\phi))$. If we restrict $\nu$ to $0 \leq \nu \leq m/2$, then $\phi_\nu - \phi \leq \pi/2$, and $\left|\sigma_\nu(x,y) \right| \geq \left| \phi_\nu - \phi \right|$, since $\cos(\sigma_\nu(x,y)) \leq \cos(\phi_\nu-\phi)$.  On the other hand, if $\phi_\nu-\phi > \pi/2$, then $\sigma_\nu(x,y) > \pi/2$.  The first inequality follows, since
\begin{align*}
	& \sum_{\nu=0}^m \frac{1}{\left| \sigma_\nu(x,y) \right| + m^{-1}}
	 \leq \sum_{\nu=0}^{\lfloor \tfrac{m}{2} \rfloor} \frac{1}{|\phi_\nu - \phi| + m^{-1}} 
	+ \frac{m}{2} \frac{2}{\pi}
	 \leq c m \log(m+1) .
\end{align*}
The proof of the second inequality is similar to the first.  Recall that $\pi - \sigma_\nu(x,y) = \sigma_\nu(-x,-y)$, and write $(-x,-y) = (r \cos( \phi) , r \sin(\phi))$, where $\phi \in ( \pi -\tfrac{\pi}{4m+2}, \pi + \tfrac{\pi}{4m+2} )$.  A similar argument shows $\sigma_\nu(-x,-y) > \pi/2$ for $\nu \leq m/2$, while $\sigma_\nu(-x,-y) > |\phi-\phi_\nu|$ for $m/2 < \nu \leq m$, and the remainder of the proof is identical to the proof of the first inequality. 
\end{proof}

\begin{prop}
	For $0 \leq \phi_1 < \phi_2 \leq \pi$,
	\begin{displaymath}
		\int_{\phi_1}^{\phi_2} 
		\frac{1}{ \left( \left\{ \left| \theta-\phi_1 \right| , \left| \phi_2 - \theta \right|\right\} + m^{-1}\right)^k} \ d\theta 
		\leq \left\{
		\begin{array}{ll}
			c \log(m+1) \ & \ k=1\\
			c m \ & \  k=2
		\end{array}
		\right.
	\end{displaymath}
\end{prop}
\begin{proof}
	The proposition follows from a change of variables in the integral.
\end{proof}

We introduce new notation to simplify the proof of the remaining estimates.  The notation $\CI_\nu(\phi_1,\phi_2)$ denotes the set of indices $\nu$ such that $\phi_1 \leq \sigma_\nu(x,y) \leq \phi_2$, and the symbol $\CI_j(\phi_1, \phi_2)$ denotes the equivalent set of indices such that $\phi_1 \leq \theta_{j} \leq \phi_2$.

Combining $A_1^+(x,y)$ and $A_1^-(x,y)$, we obtain
\begin{align*}
	&\left|A_1^+(x,y) - A_1^-(x,y) \right|\\
	& = \frac{ 4 re^{i\theta} \sin \theta_{j}  \sin \sigma_\nu(x,y) }
		  { (1 -2re^{i\theta}\cos(\theta_{j}+\sigma_\nu(x,y)) + r^2e^{2i\theta})
		    (1 -2re^{i\theta}\cos(\theta_{j}-\sigma_\nu(x,y)) + r^2 e^{2i\theta}) }.
\end{align*}
Upon substituting this into \eqref{B2Mestimate} and using \eqref{easyfactor} and \eqref{easyapprox}, we are left with approximating
\begin{align}
	&\label{A1estimate}
	\frac{1}{(2m+1)^2} \int_0^{\pi} \sum_{\nu=0}^m \sum_{j=1}^m
	\frac{1}{\sin \theta +m^{-1}} \sin^2  \theta_{j}  \\
	& \qquad \times \frac{ 1 }
	{\left( \left| \sin \left( \frac{ \theta + \theta_{j} + \sigma_\nu(x,y) }{2} \right) \right| + m^{-1} \right)
	\left( \left| \sin \left( \frac{ \theta - \theta_{j} - \sigma_\nu(x,y) }{2} \right) \right| + m^{-1} \right)} 
	\notag\\
	& \qquad \times \frac{1}
	{\left( \left| \sin \left( \frac{ \theta + \theta_{j} - \sigma_\nu(x,y) }{2} \right) \right| + m^{-1} \right)
	\left( \left| \sin \left( \frac{ \theta - \theta_{j} + \sigma_\nu(x,y) }{2} \right) \right| + m^{-1} \right) }
	\ d\theta \notag.
\end{align}
In order to approximate the sine functions, the integral over $[0,\pi]$ is divided into integrals over three subintervals: $\left[ 0, \tfrac{\pi}{4} \right]$, $\left[ \tfrac{\pi}{4} , \tfrac{3\pi}{4} \right]$, and $\left[ \tfrac{3\pi}{4}, \pi \right]$.  We will use the notation $\CH_1^1$, $\CH_1^2$, and $\CH_1^3$ to denote the expression \eqref{A1estimate} restricted over these respective sub-intervals.

\vskip .1 in
\textbf{Case 1: $0 \leq \theta \leq \pi/4$.}
With this restriction on $\theta$, the sine functions in \eqref{A1estimate} are approximated by 
\begin{align}
	& \left| \sin \left( \frac{\theta \pm_1 \theta_{j} \pm_2 \sigma_\nu(x,y)}{2} \right) \right| 
	\approx \left| \theta \pm_1 \theta_{j} \pm_2 \sigma_\nu(x,y) \right|
	\label{sinapprox2}\\
\end{align}
to obtain
\begin{align}
	\label{sineapprox1}
	\CH_1^1 & \leq c \frac{1}{(2m+1)^2} \int_0^{\tfrac{\pi}{4}}
	\sum_{\nu=0}^m \sum_{j=1}^m
	\frac{1}{\theta+m^{-1}} \theta_{j}^2\\
	& \notag \qquad \times \frac{1}
	{\left( \left| \theta + \theta_{j} + \sigma_\nu(x,y) \right| + m^{-1} \right)
	\left( \left| \theta - \theta_{j} - \sigma_\nu(x,y) \right| + m^{-1} \right)}\\
	& \notag \qquad \times \frac{1}
	{\left( \left| \theta + \theta_{j} - \sigma_\nu(x,y) \right| + m^{-1} \right)
	\left( \left| \theta - \theta_{j} + \sigma_\nu(x,y) \right| + m^{-1} \right) }.
\end{align}
First, use the inequality $\theta_{j}/(\theta + \theta_{j} + \sigma_\nu(x,y) + m^{-1}) < 1$.  Considering $\nu \in \CI_\nu \left( 0 , \theta \right)$,  we use the inequality $\theta_{j}/(\theta_{j} + \theta - \sigma_\nu(x,y) + m^{-1})<1$, and split the sum in $j$ into $\CI_j \left( 0, \theta \right)$ and $\CI_j \left( \theta , \tfrac{\pi}{2} \right)$ to obtain estimates of type 5.  For $\nu \in \CI_\nu \left( \theta , \pi \right)$, the inequality $\theta_{j}/(\theta_{j} + \sigma_\nu(x,y) - \theta + m^{-1}) <1$ is used.  Splitting the sum in $j$ into $\CI_j \left( 0 , \sigma_\nu(x,y) \right) $ and $\CI_j \left( \sigma_\nu(x,y) , \tfrac{\pi}{2} \right) $ yields two estimates of type 3.

\vskip .1 in
\noindent

\textbf{Case 2: $ \pi/4 \leq \theta \leq 3 \pi /4$.}  In this case, the factor of $( \sin(\theta) + m^{-1} )^{-1}$ in \eqref{A1estimate} is bounded away from zero and may be ignored.  We approximate one of the sine functions in \eqref{A1estimate} by
\begin{equation}
	\left| \sin \left( \frac{\theta + \theta_{j} + \sigma_\nu(x,y)}{2}\right) \right|
	\approx \left| 2\pi - \theta - \theta_{j} - \sigma_\nu(x,y) \right|
	\label{sinapprox3}
\end{equation} 
and the remaining three are estimated as in \eqref{sinapprox2} to obtain
\begin{align}
	\CH_1^2 & \leq \label{sineapprox2} 
	c \frac{1}{(2m+1)^2} \int_{\tfrac{\pi}{4}}^{\tfrac{3\pi}{4}}
	\sum_{\nu=0}^m \sum_{j=1}^m
	\theta_{j,2m}^2 \\
	& \notag \qquad \times \frac{1}
	{\left( \left| 2 \pi - \theta - \theta_{j} - \sigma_\nu(x,y) \right| + m^{-1} \right)
	\left( \left| \theta - \theta_{j} - \sigma_\nu(x,y) \right| + m^{-1} \right)} \notag\\
	& \notag \qquad \times \frac{1}
	{\left( \left| \theta + \theta_{j} - \sigma_\nu(x,y) \right| + m^{-1} \right)
	\left( \left| \theta - \theta_{j} + \sigma_\nu(x,y) \right| + m^{-1} \right) }\ d\theta.
\end{align}
First consider $\nu \in \CI_\nu \left( 0 , \theta \right)$.  Under this restriction, $2\pi - \theta -\theta_{j} - \sigma_\nu (x,y) > \frac{3 \pi}{2} - 2\theta$ and $\theta_{j}/(\theta_{j}+\theta-\sigma_\nu(x,y) +m^{-1}) < 1$.  Splitting the sum in $j$ in $\CI_j \left( 0 , \theta \right)$ and $\CI_j \left( \theta , \tfrac{\pi}{2} \right)$ yields estimates of type 5.  If $\nu \in \CI_\nu \left( \theta , \pi \right)$, first note that $\theta_{j} / (\theta_{j} + \sigma_\nu(x,y) - \theta + m^{-1}) \leq 1$, and $\theta-\theta_{j}+\sigma_\nu(x,y) \geq 2(\theta - \pi/4)$.  Substituting $2m+1-j$ for $j$ and splitting the sum in $j$ into $\CI_j \left( \frac{\pi}{2} , \theta \right) $ and $\CI_j \left( \theta , \pi \right)$ yields estimates of type 5.

\vskip .1 in
\noindent
\textbf{Case 3: $3 \pi/4 \leq \theta \leq \pi$.}  We may again ignore the factor of $\left( \sin(\theta) + m^{-1} \right)^{-1}$ in \eqref{A1estimate}.  Two of the sine functions are approximated by 
\begin{equation}
	\sin\left( \frac{\theta-\theta_{j}+\sigma_\nu(x,y)}{2} \right) \approx 2\pi-\theta+\theta_{j}-\sigma_\nu(x,y),
	\label{sinapprox6}
\end{equation}
and \eqref{sinapprox3}, and the remaining two sine functions are approximated by \eqref{sinapprox2} to obtain
\begin{align}
	\CH_1^3 & \leq \label{sineapprox3}
	c \frac{1}{(2m+1)^2} \int_{\tfrac{3\pi}{4}}^\pi 
	\sum_{\nu=0}^m \sum_{j=1}^m \theta_{j}^2  \\
	& \qquad \times \frac{1}
	{\left( \left| 2 \pi - \theta - \theta_{j} - \sigma_\nu(x,y) \right| + m^{-1} \right)
	\left( \left| \theta - \theta_{j} - \sigma_\nu(x,y) \right| + m^{-1} \right)}
	 \notag \\
	& \qquad \times \frac{1}
	{\left( \left| \theta + \theta_{j} - \sigma_\nu(x,y) \right| + m^{-1} \right)
	\left( \left| 2 \pi - \theta + \theta_{j} - \sigma_\nu(x,y) \right| + m^{-1} \right) } d\theta.
	\notag
\end{align}
First observe that $\theta + \sigma_\nu(x,y) \leq 2\pi$, and hence $\theta_{j}/ (2\pi - \theta + \theta_{j} - \sigma_\nu(x,y) + m^{-1}) < 1$.  Considering $\nu \in \CI_\nu \left( 0 , \theta \right)$, use the inequality $\theta_{j} / (\theta + \theta_{j} - \sigma_\nu(x,y) + m^{-1}) \leq 1$, and substitute $2m+1-j$ for $j$.  After splitting the sum in $j$ into $ \CI_j \left( \frac{\pi}{2} , \sigma_\nu(x,y) \right) $ and $\CI_j \left( \sigma_\nu(x,y) , \pi \right)$, estimates of type 1 are obtained.  For $\nu \in \CI_\nu \left( \theta , \pi \right) $, use the inequality $\theta_{j} / \left( \theta_{j} + \sigma_\nu(x,y) - \theta + m^{-1} \right) < 1$, and split the sum in $j$ into $\CI_j \left( \frac{\pi}{2}, \theta \right)$ and $\CI_j \left( \theta, \pi \right)$.  The resulting expressions are type 6. 
\end{proof}

\subsection{Proof of Lemma \ref{lem:upperboundA2}}

\begin{proof}
	First, we let $\zeta = re^{i\theta}$.  As in the proof of Lemma \eqref{lem:upperboundA1}, we combine terms to obtain
\begin{align*}
	& \left| A_2^+(x,y) - A_2^-(x,y) \right|= 2 \sin \theta_{j}  \sin \sigma_\nu(x,y)  
	\left|1-r^2e^{2i\theta} \right| \\ 
	& \quad \times \frac{ \left| P(\zeta,\theta_{j},\sigma_\nu(x,y))\right|}
	{ \left(1-re^{i\theta}\cos(\theta_{j} + \sigma_\nu(x,y)) + (re^{i\theta})^2\right)^2
	  \left(1-re^{i\theta}\cos(\theta_{j} - \sigma_\nu(x,y)) + (re^{i\theta})^2\right)^2},
\end{align*}
where
\begin{align}
	\label{P} & P(\zeta,\theta_{j},\sigma_\nu(x,y))  = 1 - 3 \zeta^4 
	- 2 \zeta^2 (3 - 2 \sin^2 (\theta_{j}) - 2 \sin^2(\sigma_\nu(x,y))) 
	\\
	& \quad+ 8 \zeta^3 \left(1 -2 \sin^2 ( \sigma_\nu(x,y)/2 ) 
	- 2 \sin^2 (\theta_{j}/2 ) 
	+ 4 \sin^2 ( \sigma_\nu(x,y)/2 ) \sin^2 ( \theta_{j} / 2 ) \right) \notag.
\end{align}
Substituting this into \eqref{B2Mestimate} and using \eqref{easyfactor} and \eqref{easyapprox}, it remains to estimate
\begin{align}
	&\label{A2estimate}
	\frac{1}{(2m+1)^2} \int_0^{\pi} \sum_{\nu=0}^m \sum_{j=1}^m
	\left(\sin(\pi-\theta) + m^{-1} \right) \sin^2 (\theta_{j}) \\
	& \qquad \times \frac{ \left| P(\zeta,\theta_{j}, \sigma_\nu(x,y) \right|  }
	{\left( \left| \sin \left( \frac{ \theta + \theta_{j} + \sigma_\nu(x,y) }{2} \right) \right| + m^{-1} \right)^2
	\left( \left| \sin \left( \frac{ \theta - \theta_{j} - \sigma_\nu(x,y) }{2} \right) \right| + m^{-1} \right)^2} 
	\notag\\
	& \qquad \times \frac{1}
	{\left( \left| \sin \left( \frac{ \theta + \theta_{j} - \sigma_\nu(x,y) }{2} \right) \right| + m^{-1} \right)^2
	\left( \left| \sin \left( \frac{ \theta - \theta_{j} + \sigma_\nu(x,y) }{2} \right) \right| + m^{-1} \right)^2 }
	\ d\theta \notag.
\end{align}

We will again split this integral into three sub-integrals over $\left[ 0, \tfrac{\pi}{4} \right]$, $\left[ \tfrac{\pi}{4}, \tfrac{3\pi}{4} \right]$, and $\left[ \tfrac{3\pi}{4}, \pi \right]$, and denote the part of \eqref{A2estimate} associated with these subintervals by $\CH_2^1$, $\CH_2^2$, and $\CH_2^3$, respectively.  The crucial part of the estimate is suitably approximating $\left| P(\zeta,\theta_{j}, \sigma_\nu(x,y) )\right|$.  Several different approximations will be used.  These approximations are given in the following lemma, and are referenced as needed.

\begin{lem}
	The function $P(\zeta, \theta_j, \sigma_\nu(x,y))$ satisfies the following inequalities.  
	\begin{enumerate}
	\item[(E1)]For $0 \leq \theta \leq \pi/2$,
		\begin{align}
			\label{Pest1}
			& \left| P(\zeta, \theta_{j}, \sigma_\nu(x,y)) \right| 
			\leq c \ \left( \left( \theta+\theta_{j}\ +  m^{-1} \right)^2 
			\left( \left| \theta-\theta_{j}\right| + m^{-1} \right) 
			 \right. \\
			& \qquad \left. \qquad + (\sigma_\nu(x,y))^2(\theta + (\sigma_\nu(x,y))^2 + \theta_{j}^2 + m^{-1}) \right). \notag
		\end{align}
	\item[(E2)]  For $0 \leq \theta \leq \pi/2$,
		\begin{align}
			\label{Pest2}
			& \left| P(\zeta, \theta_{j}, \sigma_\nu(x,y) \right| \leq c \ ( (\theta + m^{-1})^3 + (\theta_{j} + \sigma_\nu(x,y))^2(|\theta_{j}-\sigma_\nu(x,y)|)^2  \\
			& \qquad \qquad +(\theta+m^{-1})(\theta_{j}^2 +(\sigma_\nu(x,y))^2)) \notag.
		\end{align}
	\item[(E3)] For $\pi/4 \leq \theta \leq 3\pi/4$,
		\begin{align}
			& \label{Pest4}
			\left| P(\zeta, \theta_{j}, \sigma_\nu(x,y)) \right| \leq c \ \left(  \left( \left| \theta + \theta_{j}\right| + m^{-1} \right)
			\left( \left| \theta - \theta_{j}\right| + m^{-1} \right) 
			+ (\sigma_\nu(x,y))^2 \right).
		\end{align}
	\item[(E4)] For $\pi/4 \leq \theta \leq \pi$,
		\begin{align}
			\label{Pest3}
			& \left| P(\zeta, \theta_{j}, \sigma_\nu(x,y)) \right| \\ 
			&\qquad \leq c \ \left(
			\left( \left| \pi-\theta+\theta_{j} \right| + m^{-1}\right)^2
			\left( \left| \pi -\theta - \theta_{j} \right| + m^{-1}\right)
			+ (\pi-\sigma_\nu(x,y))^2 \right) \notag.
		\end{align}
	\end{enumerate}
\end{lem}
\begin{proof}
	We first prove the estimate \eqref{Pest1}.  We define
\begin{align}
	& P_1(\zeta,\theta_{j,2m}) := 1 - 2 \zeta^2 ( 3 - 2\sin^2 \theta_{j} ) 
	+ 8 \zeta^3 \cos \theta_{j}  - 3\zeta^4,\\
	& \notag P_2^\pm(\zeta,\theta_{j}, \sigma_\nu(x,y)) := 4 \zeta^2 \sin^2 \sigma_\nu(x,y) 
	\pm 16 \zeta^3 \sin^2 \tfrac{\sigma_\nu(x,y)}{2} \cos \theta_{j}  , \notag
\end{align} 
so we may write $P(\zeta,\theta_{j},\sigma_\nu(x,y)) = P_1(\zeta, \theta_{j}) + P_2^-(\zeta, \theta_{j}, \sigma_\nu(x,y))$.  It is possible to factor $P_1(\zeta,\theta_{j})$ as
\begin{equation}
	\label{P1factor}
	P_1(\zeta,\theta_{j}) = (1+2 \zeta \cos \theta_{j} - 3\zeta^2)(1-2\zeta\cos \theta_j + \zeta^2).
\end{equation}
The second factor of \eqref{P1factor} is approximated by \eqref{easyfactor} and \eqref{easyapprox} as before, and the first factor may be further factored as
\begin{align}
	& \label{P1factor2}
	(1 + 2 \zeta \cos \theta_{j} - 3 \zeta^2) = -3 \left(\zeta + \frac{1}{3} \left(\sqrt{4 -\sin^2 \theta_{j} }-\cos \theta_{j}  \right)\right)\\
	&  \notag \qquad \qquad \times
	     \left( \zeta - \frac{1}{3} \left(\sqrt{4 -\sin^2 \theta_{j} }+\cos \theta_{j}  \right)\right).
\end{align}
The first factor of \eqref{P1factor2} will not be used.  Using the double angle identity for cosine, the second factor of \eqref{P1factor2} is approximated by
\begin{align*}
	& \left| \zeta - \frac{1}{3} \left(\sqrt{4 -\sin^2 \theta_{j} }+\cos \theta_{j} \right) \right|\\
	& \leq c \left( 3 \sin \theta  + 
	\left| 3 \cos\theta - \cos \theta_{j}  -\sqrt{ 4 - \sin^2 \theta_{j} } \right| + m^{-1} \right)\\
	& \leq c \left( \theta + \theta_{j}^2 + \left| 2 - \sqrt{4-\sin^2 \theta_{j} } \right| +m^{-1}\right).
\end{align*}
Since $2 - \sqrt{4-\sin^2 \theta_{j} } \leq \sin \theta_{j} $, we obtain 
\begin{equation}
	\label{P1one}
	P_1(\zeta, \theta_{j}) \leq c \left( \left(\left| \frac{\theta+\theta_{j}}{2}\right| + m^{-1}\right) 
	\left( \left| \frac{\theta-\theta_{j}}{2} \right| + m^{-1}\right)  
	\left(\theta+\theta_{j} + m^{-1}\right) \right).  
\end{equation}
Now considering $P_2^-(\zeta,\theta_{j})$, the double angle identities for sine and cosine are used to obtain
\begin{align}
	\label{P2one}
	\left| P_2^-(\zeta, \theta_{j}, \sigma_\nu(x,y) \right| & \leq c \left(\sigma_\nu(x,y) \right)^2  \left| \cos^2 \tfrac{\sigma_\nu(x,y)}{2} -  z \cos \theta_{j} \right|\\
	\notag & \leq c (\sigma_\nu(x,y))^2 ( \theta + (\sigma_\nu(x,y))^2 + \theta_{j}^2 + m^{-1}).
\end{align}
Adding the estimates \eqref{P1one} and \eqref{P2one} , we arrive at the estimate \eqref{Pest1}.

The proof of the estimate \eqref{Pest4} follows from replacing $\theta$ with a constant in \eqref{Pest1}.

We next prove the estimate \eqref{Pest2}.  We first re-write $P(\zeta,\theta_{j}, \sigma_\nu(x,y))$ as
\begin{align*}
	& P(\zeta,\theta_{j}, \sigma_\nu(x,y)) = 1 - 3\zeta^4 - 6 \zeta^2 +8 \zeta^3
	+ 4\zeta^2 \bigg[ \sin^2 \theta_{j} + \sin^2 \sigma_\nu(x,y) \\ 
	& \qquad  + 4 \left( 2 \sin^2 \tfrac{\theta_{j}}{2}  \sin^2 \tfrac{\sigma_\nu(x,y)}{2} 
	-\sin^2 \tfrac{\theta_{j}}{2} -\sin^2 \tfrac{\sigma_\nu(x,y)}{2} \right)\\
	& \qquad 
	+ 4(\zeta-1) \left( 2 \sin^2 \tfrac{\theta_{j}}{2} \sin^2 \tfrac{\sigma_\nu(x,y)}{2}
	-\sin^2 \tfrac{\theta_{j}}{2} - \sin^2 \tfrac{\sigma_\nu(x,y)}{2}  \right) \bigg].  
\end{align*}
It is easily checked that $1-3\zeta^4-6\zeta^2+8\zeta^3 = -(\zeta-1)^3(3\zeta+1)$, and applying the double angle identity for sines,
\begin{align*}
	& \sin^2 \theta_{j} + \sin^2 \sigma_\nu(x,y) 
	+ 4 \left( 2 \sin^2 \tfrac{\theta_{j}}{2}  \sin^2 \tfrac{\sigma_\nu(x,y)}{2}  
	- \sin^2 \tfrac{\theta_{j}}{2} -\sin^2 \tfrac{\sigma_\nu(x,y)}{2} \right)\\ 
	& = -4 \left( \sin^2 \tfrac{\theta_{j}}{2} - \sin^2 \tfrac{\sigma_\nu(x,y)}{2} \right)^2.
\end{align*}
Finally, we may approximate the sine functions to obtain \eqref{Pest2}.

Finally, we prove the estimate \eqref{Pest3}.  First, replace $\sin^2 \sigma_\nu(x,y) $ with $\sin^2 (\pi-\sigma_\nu(x,y))$ and $\sin^2 \tfrac{\sigma_\nu(x,y)}{2}$ with $1 - \sin^2 \tfrac{\pi-\sigma_\nu(x,y)}{2}$ in \eqref{P} to obtain
\begin{align}
	\label{Pfactor1} \left| P(\zeta, \theta_{j}, \sigma_\nu(x,y)) \right| \leq & \left| 1 - 3\zeta^4 -4\zeta^2 (3 - 2 \sin^2 \theta_{j}) + 8 \zeta^3(2 \sin^2 \tfrac{\theta_{j}}{2} -1)  \right| \\
	& + c \left| P_2^+(\zeta, \theta_{j}, \pi-\sigma_\nu(x,y)) \right| \notag.
\end{align}
The inequality $|P_2^+(\zeta,\theta_{j}, \pi-\sigma_\nu(x,y))| < c (\pi - \sigma_\nu(x,y))^2$ follows easily from the definition of $P_2^+$.  The first term in \eqref{Pfactor1} becomes $|P_1(\zeta, \pi-\theta_{j})|$, after replacing $3 - 2 \sin^2 \theta_{j} $ with $3 - 2 \sin^2(\pi-\theta_{j})$ and $2 \sin^2 \tfrac{\theta_{j}}{2} -1 $ with $1 - 2 \sin^2 \left( \tfrac{\pi-\theta_{j}}{2}\right)$.  $P_1(\zeta,\pi-\theta_{j})$ factors as in \eqref{P1factor}.  The factor of $1+2\cos(\theta_{j})\zeta -3\zeta^2$ can be factored further, as 
\begin{align*}
	&1+2\zeta\cos(\pi-\theta_{j})-3\zeta^2 = 1-2\zeta\cos \theta_{j} -3\zeta^2 \\
	&= -3\left(\zeta+\frac{1}{3}\left( \cos\theta_{j} - \sqrt{4-\sin^2\theta_{j}}\right) \right)\left(\zeta+\frac{1}{3}\left( \cos \theta_{j} +\sqrt{4-\sin^2 \theta_{j} }\right) \right).
\end{align*}
The first factor will not be used, but the absolute value of the second factor may be approximated by
\begin{align*}
	& \left| \zeta + \frac{1}{3} \left( \cos \theta_{j} + \sqrt{4-\sin^2 \theta_{j}} \right) \right|\\
	& \leq c \left( |\sin \theta | + \left| \cos \theta  + \frac{1}{3} \left( \cos \theta_{j} + \sqrt{4-\sin^2 \theta_{j}} \right) \right| + m^{-1}\right)\\
	& \leq c\left( \pi-\theta + \left| 6 \sin^2 \tfrac{\pi-\theta}{2} - 2 \sin^2 \tfrac{\theta_{j}}{2} + \sqrt{4 -\sin^2 \theta_{j} } -2\right| + m^{-1} \right)\\
	& \leq c \left( \pi-\theta + \theta_{j} + m^{-1} \right),
\end{align*}
and we are able to obtain the estimate \eqref{Pest3}.
\end{proof}

\vskip .1 in
\noindent
\textbf{Case 1: $0 \leq \theta \leq \pi/4$.}  
Approximating the sine functions in \eqref{A2estimate} with \eqref{sinapprox2}, we obtain
\begin{align*}
	\CH_2^1 \leq & \ c \frac{1}{(2m+1)^2} \int_0^{\tfrac{\pi}{4}}
	\sum_{\nu=0}^m \sum_{j=1}^m
	 \theta_{j,2m}^2\\
	& \times \frac{ \left| P(\zeta, \theta_{j}, \sigma_\nu(x,y) \right| }
	{\left( \left| \theta + \theta_{j} + \sigma_\nu(x,y) \right| + m^{-1} \right)^2
	\left( \left| \theta - \theta_{j} - \sigma_\nu(x,y) \right| + m^{-1} \right)^2}\\
	& \times \frac{1}
	{\left( \left| \theta + \theta_{j} - \sigma_\nu(x,y) \right| + m^{-1} \right)^2
	\left( \left| \theta - \theta_{j} + \sigma_\nu(x,y) \right| + m^{-1} \right)^2 } \ d\theta.
\end{align*}
First considering $\nu \in \CI_\nu \left( 0 , \theta \right) $, note that $\theta_{j}^2 / (\theta + \theta_{j} - \sigma_\nu(x,y) + m^{-1})^2 <1$.  Approximate $\left| P(\zeta, \theta_{j}, \sigma_\nu(x,y) \right|$ using \eqref{Pest1}.  For the first term in the sum in \eqref{Pest1}, use the inequality $( \theta+\theta_{j} + m^{-1})^2 /(\theta+\theta_{j} + \sigma_\nu(x,y) + m^{-1})^2 < 1$, and then split the sum in $j$ into $\CI_j \left( 0 , \theta \right)$ and $\CI_j \left( \theta , \tfrac{\pi}{2} \right) $ to obtain two estimates of type 7.  For the second term in \eqref{Pest1}, use the inequality $ (\theta + (\sigma_\nu(x,y))^2 + \theta_{j}^2 + m^{-1})/ (\theta + \theta_{j} + \sigma_\nu(x,y) + m^{-1})^2 < c (\theta+m^{-1})^{-1}$, and split the sum in $j$ into $\CI_j \left( 0 , \theta \right) $ and $\CI_j \left( \theta , \tfrac{\pi}{2} \right) $ to obtain two estimates type 4.

For $ \nu \in \CI_\nu \left( \theta , \pi \right) $, use the inequality $ \theta_{j}^2 / ( \theta_{j} + \sigma_\nu(x,y) - \theta + m^{-1})^2 < 1$, and approximate $ | P(\zeta,\theta_{j}, \sigma_\nu(x,y) ) |$ using \eqref{Pest2}.  For the first term of the sum in \eqref{Pest2}, use the inequality $(\theta+m^{-1})^2/(\theta+\theta_{j} + \sigma_\nu(x,y) +m^{-1})^2 < 1$ and then split the sum in $j$ into $\CI_j \left( 0 , \sigma_\nu(x,y) \right) $ and $\CI_j \left( \sigma_\nu(x,y) , \tfrac{\pi}{2} \right) $ to obtain two estimates of type 4.  For the second term of the sum in \eqref{Pest2}, first use the fact that $(\theta_{j} + \sigma_\nu(x,y))^2 / ( \theta+\theta_{j}+\sigma_\nu(x,y)+m^{-1})^2 < 1$, and then split the sum in $j$ into $\CI_j \left( 0 , \sigma_\nu(x,y) \right) $ and $\CI_j \left( \sigma_\nu(x,y) , \tfrac{\pi}{2} \right) $ to obtain two estimates of type 2.  Finally, for the third term in \eqref{Pest2}, use the inequality $( (\sigma_\nu(x,y))^2 + \theta_{j}^2) / (\theta+ \theta_{j} + \sigma_\nu(x,y) + m^{-1})^2 < 1$, and then split the sum in $j$ into $\CI_j \left( 0 , \sigma_\nu(x,y) \right) $ and $\CI_j \left( \sigma_\nu(x,y) , \tfrac{\pi}{2} \right) $ to obtain two estimates of type 4.  

\vskip .1 in
\noindent
\textbf{Case 2: $\pi/4 \leq \theta \leq 3 \pi /4$:}  If $\sigma_\nu(x,y) \leq \pi/2$, the sine functions in \eqref{A2estimate} may be approximated as in Case 1, while if $\sigma_\nu(x,y) \geq \pi/2$, one of the sine functions may be approximated by \eqref{sinapprox3} and the remaining three by \eqref{sinapprox2}, to obtain 
\begin{align*}
	\CH_2^2 \leq & \ c  \frac{1}{(2m+1)^2} \int_{\tfrac{\pi}{4}}^{\tfrac{3\pi}{4}}
	\bigg( \sum_{\nu \in \CI_\nu \left( 0 , \tfrac{\pi}{2} \right) } \sum_{j=1}^m
	 \theta_{j}^2\\
	& \times \frac{ \left| P(\zeta, \theta_{j}, \sigma_\nu(x,y) \right| }
	{\left( \left| \theta + \theta_{j} + \sigma_\nu(x,y) \right| + m^{-1} \right)^2
	\left( \left| \theta - \theta_{j} - \sigma_\nu(x,y) \right| + m^{-1} \right)^2}\\
	& \times \frac{1}
	{\left( \left| \theta + \theta_{j} - \sigma_\nu(x,y) \right| + m^{-1} \right)^2
	\left( \left| \theta - \theta_{j} + \sigma_\nu(x,y) \right| + m^{-1} \right)^2 } \bigg)\\ 
	& + 
	\bigg(\sum_{\nu \in \CI_\nu \left( \tfrac{\pi}{2}, \pi \right) } \sum_{j=1}^m
	 \theta_{j}^2\\
	& \times \frac{ \left| P(z, \theta_{j}, \sigma_\nu(x,y) \right| }
	{\left( \left| 2\pi- \theta - \theta_{j} - \sigma_\nu(x,y) \right| + m^{-1} \right)^2
	\left( \left| \theta - \theta_{j} - \sigma_\nu(x,y) \right| + m^{-1} \right)^2}\\
	& \times \frac{1}
	{\left( \left| \theta + \theta_{j} - \sigma_\nu(x,y) \right| + m^{-1} \right)^2
	\left( \left| \theta - \theta_{j} + \sigma_\nu(x,y) \right| + m^{-1} \right)^2 } \bigg)\ d\theta.
\end{align*}

We consider first $\nu \in \CI_\nu \left( 0 , \tfrac{\pi}{2} \right) \cap \CI_\nu \left( 0 , \theta \right) $.  Under this restriction, $\theta + \theta_{j} + \sigma_\nu(x,y)  \geq \pi/4 $ and so this factor may be ignored.  Use the inequality $\theta_{j}^2 / ( \theta + \theta_{j} - \sigma_\nu(x,y) + m^{-1})^2 < 1$, and approximate $P(\zeta,\theta_{j}, \sigma_\nu(x,y))$ by \eqref{Pest4}.  After splitting the sum in $j$ into $ \CI_j \left( 0 , \theta \right) $ and $\CI_j \left( \theta , \pi \right) $, estimates of types 2 and 7 are obtained. 

Next, we consider $\nu \in \CI_\nu \left( \theta , \tfrac{\pi}{2} \right)$.  The factor of $\theta+\theta_{j} + \sigma_\nu(x,y) +m^{-1}$ may again be ignored, and the inequality $\theta_{j}^2 / ( \theta_{j}^2 + \sigma_\nu(x,y) - \theta + m^{-1})^2 < 1$ is used.  We only need to estimate $P(\zeta,\theta_{j}, \sigma_\nu(x,y))$ by a constant.  After splitting the sum in $j$ into $\CI_j \left( 0, \sigma_\nu(x,y) \right)$ and $\CI_j \left( \sigma_\nu(x,y) , \tfrac{\pi}{2} \right)$, estimates of type 2 are obtained.

Now we consider $\nu \in \CI_\nu \left( \tfrac{\pi}{2} , \theta \right) $.  For this range of $\nu$, $\theta - \theta_{j} + \sigma_\nu(x,y) \geq \pi/4$ and $\theta_{j}^2 / \left( \theta + \theta_{j}-\sigma_\nu(x,y) + m^{-1} \right)^2 < 1$.  We again estimate $|P(\zeta,\theta_{j}, \sigma_\nu(x,y)|$ by a constant.  After substituting $2m+1-j$ for $j$, we split the sum in $j$ into $\CI_j \left( \tfrac{\pi}{2} , \sigma_\nu(x,y) \right) $ and $\CI_j \left( \sigma_\nu(x,y) , \pi \right)$ to obtain estimates of type 2.

To complete this case, we consider $\nu \in \CI_\nu \left( \theta , \pi \right) \cap \CI_\nu \left( \tfrac{\pi}{2} , \pi \right)$.  Under these restrictions, $\theta - \theta_{j} + \sigma_\nu(x,y) + m^{-1} > \theta+m^{-1}$ and $\theta_{j}^2 / ( \theta_{j}+\sigma_\nu(x,y)-\theta +m^{-1} ) ^2 < 1$.  Use the estimate for $\left|P(\zeta, \theta_{j}, \sigma_\nu(x,y) \right|$ in \eqref{Pest3}, and substitute  $2m+1-j$ for $j$.  Upon splitting the sum in $j$ into $\CI_j \left( \tfrac{\pi}{2} , \theta \right) $ and  $\CI_j \left( \theta , \pi \right) $, we are left with estimates of types 4 and 2.

\vskip .1 in
\noindent
\textbf{Case 3: $3 \pi /4 \leq \theta \leq \pi$.}  Recalling \eqref{A2}, a factor of $\pi-\theta+m^{-1}$ is present in the numerator.  After approximating two of the sine functions in \eqref{A2estimate} by \eqref{sinapprox3}, and \eqref{sinapprox6}, and the remaining two sine functions by \eqref{sinapprox2}, we obtain
\begin{align*}
	\CH_2^3 \leq & \ c  \frac{1}{(2m+1)^2} \int_{\tfrac{3\pi}{4}}^{\pi}(\pi-\theta + m^{-1})
	\sum_{\nu=0}^m \sum_{j=1}^m
	 \theta_{j}^2\\
	& \times \frac{ \left| P(\zeta, \theta_{j}, \sigma_\nu(x,y) \right| }
	{\left( \left|2\pi- \theta - \theta_{j} - \sigma_\nu(x,y) \right| + m^{-1} \right)^2
	\left( \left| \theta - \theta_{j} - \sigma_\nu(x,y) \right| + m^{-1} \right)^2}\\
	& \times \frac{1}
	{\left( \left| \theta + \theta_{j} - \sigma_\nu(x,y) \right| + m^{-1} \right)^2
	\left( \left| 2\pi -\theta + \theta_{j} - \sigma_\nu(x,y) \right| + m^{-1} \right)^2 } .\\ 
\end{align*}
First, consider $\nu \in \CI_\nu \left( 0 , \tfrac{\pi}{2} \right)$.  Note that $\left| 2\pi - \theta  - \sigma_\nu(x,y) + \theta_{j} \right| \geq \pi/2$, and also $\left| \theta-\sigma_\nu(x,y) +\theta_{j} \right| > \pi/2$.  Also, note that $\left| 2\pi - \theta - \theta_{j} - \sigma_\nu(x,y) \right| > \left| \pi - \theta \right| $.  Approximating $|P(\zeta,\theta_{j}, \sigma_\nu(x,y))|$ by a constant, estimates of type 4 are obtained.

We next consider $\nu \in \CI_\nu \left( \tfrac{\pi}{2} , \theta \right)$.  We estimate $\left|P(\zeta, \theta, \theta_{j}, \sigma_\nu(x,y) ) \right|$ with \eqref{Pest3}, and use the inequality $\theta_{j}^2/ \left( \theta+ \theta_{j} - \sigma_\nu(x,y) + m^{-1} \right)^2 < 1$.  For the first term in the sum in \eqref{Pest3}, use the inequality $(\pi+\theta_{j} - \theta + m^{-1})^2 / ( 2\pi-\theta-\sigma_\nu(x,y) + \theta_{j} + m^{-1} )^2 < 1$ and split the sum in $j$ into $\CI_j \left( 0, \pi-\sigma_\nu(x,y) \right) $ and $\CI_\nu \left( \pi-\sigma_\nu(x,y), \tfrac{\pi}{2} \right)$ to obtain two estimates of type 4.  For the second term of the sum in the \eqref{Pest3}, use the inequality $(\pi-\sigma_\nu(x,y))^2 / ( 2\pi - \theta - \sigma_\nu(x,y) + \theta_{j} + m^{-1})^2 < 1$ and estimates of type 4 follow by the same splitting of the sum in $j$.

Finally, we consider $\nu \in \CI_\nu \left( \theta , \pi \right)$.  Use the inequality $\theta_{j}^2 / (\theta_{j} + \sigma_\nu(x,y)-\theta + m^{-1})^2 < 1$ and approximate $|P(\zeta, \theta_{j}, \sigma_\nu(x,y))|$ using \eqref{Pest3}.  For the first sum in \eqref{Pest3}, use the inequality $(\pi+\theta_{j}-\theta + m^{-1})^2/(2\pi - \theta - \sigma_\nu(x,y) + \theta_{j} + m^{-1} )^2 < 1$, and then split the sum in $j$ into $\CI_j \left( 0, \pi - \theta \right)$ and $\CI_j \left( \pi-\theta, \tfrac{\pi}{2} \right)$.  Utilizing the factor of $|\pi-\theta-\theta_{j}| + m^{-1}$, two estimates of type 7 are obtained.  For the second term in \eqref{Pest3}, use the inequality $( \pi - \sigma_\nu(x,y))(\pi-\theta + m^{-1}) / (2\pi - \theta - \sigma_\nu(x,y) + \theta_{j} + m^{-1})^2 < 1$ and then split the sum in $j$ into  $\CI_j \left( 0, \pi - \theta \right)$ and $\CI_j \left( \pi-\theta, \tfrac{\pi}{2} \right)$ to obtain two estimates of type 7.

\end{proof}

\end{document}